\documentclass[11pt,a4paper]{article}

\usepackage[utf8]{inputenc}
\usepackage[T1]{fontenc}
\usepackage{lmodern}
\usepackage[margin=1in]{geometry}
\usepackage{amsmath,amsfonts,amssymb,amsthm}
\usepackage{graphicx}
\usepackage{cite}
\usepackage{xcolor}
\usepackage[hidelinks]{hyperref}
\allowdisplaybreaks[1]
\numberwithin{equation}{section}

\newtheorem{theorem}{Theorem}[section]

\newtheorem{lemma}[theorem]{Lemma}
\theoremstyle{remark}
\newtheorem{remark}[theorem]{Remark}
\theoremstyle{plain}

\usepackage{multirow}
\usepackage[T1]{fontenc}
\usepackage{lmodern}
\usepackage{amsmath,amssymb,mathtools,bm,mathrsfs}
\usepackage{graphicx,algorithm,algpseudocode,booktabs,array,tabularx,microtype}
\usepackage[hidelinks]{hyperref}
\usepackage{placeins}
\usepackage{xcolor}
\usepackage{cite}

\hypersetup{colorlinks=true,linkcolor=blue,citecolor=blue,urlcolor=blue}
\newcolumntype{C}[1]{>{\centering\arraybackslash}m{#1}}

\allowdisplaybreaks[2]
\newcommand{\R}{\mathbb R}

\newcommand{\ip}[2]{\left\langle #1,#2\right\rangle}
\newcommand{\norm}[1]{\left\|#1\right\|}
\newcommand{\pos}[1]{\left[#1\right]_+}
\newcommand{\cL}{\mathcal L}

\newcommand{\proj}{\operatorname{Proj}}
\newcommand{\dom}{\operatorname{dom}}
\newcolumntype{Y}{>{\raggedright\arraybackslash}X}
\newcolumntype{P}[1]{>{\raggedright\arraybackslash}p{#1}}

\title{Augmented Lagrangian methods with Nesterov acceleration for convex optimization with nonlinear inequality constraints\thanks{This work was supported by the National Natural Science Foundation of China (Grant No. 12601606), Sichuan Science and Technology Program (Grant No. 2025ZNSFSC0813), and the Talent Introduction Project of Xihua University (Grant No. Z241102).}}

\author{Xin He\thanks{School of Science, Xihua University, Chengdu, Sichuan, P.R. China (\href{mailto:hexinuser@163.com}{\texttt{hexinuser@163.com}}).}}

\date{}
\hypersetup{
  pdftitle={Augmented Lagrangian methods with Nesterov acceleration for convex optimization with nonlinear inequality constraints},
  pdfauthor={Xin He}
}

\begin{document}
\maketitle

\begin{abstract}
In this paper, we develop Nesterov accelerated augmented Lagrangian methods for composite convex optimization with nonlinear inequality constraints. The methods are derived through a compatible discretization of inertial primal-dual dynamics, yielding a unified framework that couples difference extrapolation of the constraint mapping with a projected multiplier. Under convexity and weighted summability of the primal stationarity errors, we prove boundedness of the primal-dual iterates and simultaneous $\mathcal O(k^{-2})$ last-iterate rates for the objective residual and constraint violation. Under strict parameter conditions, the primal-dual sequence converges to a Karush-Kuhn-Tucker pair, and both residuals improve to $o(k^{-2})$. These results require no strong convexity of the objective or constraints. Numerical experiments examine the performance of the proposed methods.
\end{abstract}



\section{Introduction}
\subsection{Problem formulation}
 We consider the composite convex problem with nonlinear inequality constraints:
\begin{equation}\label{ques_main}
 \min_{x\in\R^n} f(x):=\phi(x)+h(x)
 \qquad\text{subject to}\qquad g(x)\leq0,
\end{equation}
where $\phi:\R^n\to\R$ is convex with $L_\phi$-Lipschitz continuous gradient, $h:\R^n\to(-\infty,+\infty]$ is proper, closed, and convex, and
\[
 g=(g_1,\ldots,g_m)^\top:\R^n\to\R^m,
 \quad\text{with}\quad g_i:\R^n\to\R,\quad i=1,\ldots,m.
\]
Each $g_i$ is convex and continuously differentiable. Vector inequalities and positive parts are componentwise. All vector norms are Euclidean, and matrix norms are induced norms unless specified otherwise. This formulation provides a common framework for a range of constrained convex optimization models. Examples include convex quadratically constrained quadratic programming \cite{XuMathProg,XuLinear2021}, constrained sparse recovery \cite{Berg2011}, and Neyman-Pearson classification with smooth convex surrogate losses \cite{LinNP2018,Rigollet}.

Define the Lagrangian function associated with \eqref{ques_main} by
\[
 \cL(x,\lambda):=f(x)+\ip{\lambda}{g(x)},
 \qquad (x,\lambda)\in\R^n\times\R_+^m.
\]
Denote the Jacobian of the constraint mapping by
\[J_g(x):=[\nabla g_1(x)\ \cdots\ \nabla g_m(x)]^\top
\in\R^{m\times n}.\]
Let $\Omega$ denote the set of pairs
$(x^\star,\lambda^\star)\in\dom h\times\R^m$
satisfying the Karush-Kuhn-Tucker (KKT) conditions
\begin{equation}\label{eq:KKT}
\left\{
\begin{aligned}
 &0\in\nabla\phi(x^\star)+\partial h(x^\star)
       +J_g(x^\star)^\top\lambda^\star,\\
 &g_i(x^\star)\leq0,\qquad \lambda_i^\star\geq0,
       &&i=1,\ldots,m,\\
 &\lambda_i^\star g_i(x^\star)=0,
       &&i=1,\ldots,m,
\end{aligned}
\right.
\end{equation}
where $\partial h$ is the convex subdifferential of $h$.
Throughout, we assume that $\Omega\neq\varnothing$ and denote
the optimal value by $f^\star$. This assumption holds, for example,
if a minimizer exists and Slater's condition is satisfied:
there is a point $\widetilde x\in\operatorname{ri}(\dom h)$
such that $g_i(\widetilde x)<0$ for all $i=1,\ldots,m$
\cite{Rockafellar,Bauschke16}. Under convexity, each $(x^\star,\lambda^\star)\in\Omega$
consists of a primal optimal solution and an associated optimal
multiplier, and satisfies the saddle-point relation
\begin{equation}\label{eq:saddle}
 \cL(x^\star,\lambda)
 \leq\cL(x^\star,\lambda^\star)
 =f^\star
 \leq\cL(x,\lambda^\star),\quad \forall (x,\lambda)\in\R^n\times\R_+^m.
\end{equation}

Nesterov acceleration gives $\mathcal O(k^{-2})$ objective rates for smooth and composite unconstrained convex optimization \cite{NesterovSovMath,NesterovBook,BeckSiis,HeCoap}. Suitable extrapolation gives iterate convergence \cite{ChambolleJota,AttouchMathProg} and can improve the objective rate to $o(k^{-2})$ \cite{AttouchPtSiopt}. Accelerated methods have also been developed for problems with affine equality constraints \cite{HeAuto,BotMathProg,HeArxiv,HeNA,LuoMc,XuSIOPT}. For problem \eqref{ques_main} with nonlinear constraints, existing accelerated methods differ in their smoothness assumptions, in whether they require strong convexity of the objective or the constraints, and in the primal sequences for which rates are established, such as generated iterates and weighted primal averages \cite{HamedaniAybat2021,ConEx2023,ZhuSiopt,LinNeurips,ZhangLanACGD,DengInformsJoc}. For problem~\eqref{ques_main}, we seek methods with Nesterov acceleration that provide accelerated objective and feasibility estimates at the generated primal iterates, together with convergence of the primal-dual sequence. This leads to the following question.
\begin{quote}
\emph{Can a class of augmented Lagrangian methods with Nesterov acceleration for problem~\eqref{ques_main} achieve, under different parameter conditions, simultaneous $\mathcal O(k^{-2})$ last-iterate rates for the objective residual and constraint violation, and $o(k^{-2})$ rates for both residuals together with convergence of the  primal-dual sequence, without assuming strong convexity of the objective or the constraints?}
\end{quote}

To address this question, we discretize the continuous-time primal-dual dynamical system in \cite{HeNonlinear2026} and obtain a class of inexact  accelerated augmented Lagrangian methods. The methods extrapolate the constraint values and use a projected multiplier. A discrete energy estimate then yields accelerated residual rates and primal-dual convergence under explicit parameter and error conditions.

\subsection{Related work}

\paragraph{Nesterov acceleration.} Nesterov's accelerated gradient method
\cite{NesterovSovMath,NesterovBook} combines extrapolation with
gradient descent to accelerate smooth unconstrained convex
minimization. The fast iterative shrinkage-thresholding algorithm
(FISTA) \cite{BeckSiis} extends this approach to composite problems
of the form
\[ \min_{x\in\R^n} \phi(x)+h(x),\]
where $\phi:\R^n\to\R$ is convex with an
$L_\phi$-Lipschitz continuous gradient, and
$h:\R^n\to(-\infty,+\infty]$ is proper, closed, and convex.
Assume that the solution set is nonempty.
A FISTA-type iteration with a modified extrapolation rule,
as studied in \cite{ChambolleJota,AttouchMathProg,SuJmlr}, takes the form
\begin{equation}\label{eq:fista}
\begin{aligned}
 y_k&=x_k+\frac{k-1}{k+\alpha-1}(x_k-x_{k-1}),\\
 x_{k+1}&=\operatorname{prox}_{\tau h}
 \bigl(y_k-\tau\nabla\phi(y_k)\bigr),
\end{aligned}
\end{equation}
where $\alpha\geq3$, $0<\tau\leq1/L_\phi$, and
$\operatorname{prox}_{\tau h}(v)
=\arg\min_x\{h(x)+\norm{x-v}^2/(2\tau)\}$.
For $\alpha\geq3$, this iteration achieves an
$\mathcal O(k^{-2})$ objective rate \cite{AttouchMathProg,SuJmlr}.
For $\alpha>3$, the rate improves to $o(k^{-2})$ \cite{AttouchPtSiopt}, and the iterates
converge to a minimizer \cite{ChambolleJota,AttouchMathProg}.

 For affine equality constraints, Lagrangian-based frameworks provide  objective  and feasibility estimates, including accelerated estimates for strongly convex objectives \cite{Sabach2022,LuoJogo,TaoMC}. Without strong convexity, accelerated primal-dual methods attain simultaneous $\mathcal O(k^{-2})$ rates \cite{HeAuto,HeNA,XuSIOPT,LuoMc}. Accelerated augmented Lagrangian methods combine such rates with primal-dual convergence under suitable inertial rules \cite{BotMathProg}, while implicit and partially explicit variants also admit $o(k^{-2})$ residual refinements and primal-dual convergence \cite{HeArxiv}. Extending this combination of rates and iterate convergence to the nonlinear inequality constrained problem \eqref{ques_main} is the focus of the present work.
 
\paragraph{Augmented Lagrangian methods for nonlinear inequality constraints.}
Augmented Lagrangian methods for problem~\eqref{ques_main} originate in classical multiplier methods \cite{Powell1969,Hestenes1969,RockALM1976} and connect naturally with convex-concave saddle-point optimization \cite{Nemirovski2004,ChambollePock2011,Chambolle}. Xu \cite{XuLinear2021} establishes primal-dual convergence and $\mathcal O(k^{-1})$ residual estimates for averaged primal iterates of a linearized augmented Lagrangian method. Hamedani and Aybat \cite{HamedaniAybat2021} treat nonlinear convex-concave couplings, obtaining an ergodic gap rate and iterate convergence in the convex regime. For inexact augmented Lagrangian methods, Xu \cite{XuMathProg} and Lu and Zhou \cite{LuSiopt} analyze the first-order complexity, including the cost of primal subproblem solves. Related developments include proximal augmented Lagrangian methods for nonconvex objectives with convex constraints \cite{KongMOR2023,DPALM2026}, conically constrained variational inequalities \cite{ALAVI2025}, conditional-gradient methods \cite{LanSiopt}, finite-sum level-set methods \cite{LinNP2018}, and level-constrained proximal methods \cite{BoobMathProgb}.
 
Strong convexity yields faster guarantees in several primal-dual frameworks. In the smooth  setting, Boob et al. \cite{ConEx2023} gives $\mathcal O(k^{-2})$ objective and feasibility estimates for averaged primal iterates when the  objective function is strongly convex. Hamedani and Aybat \cite{HamedaniAybat2021} obtain an accelerated ergodic gap rate when the primal objective component is strongly convex and the coupling is linear in the multiplier. Zhu et al. \cite{ZhuSiopt} use a generalized augmented Lagrangian to obtain $\mathcal O(k^{-1})$ primal residual bounds for averaging and last-iterate variants, with $\mathcal O(k^{-2})$ rates under strong convexity of  objective function and the corresponding regularity conditions. Lin and Deng \cite{LinNeurips} instead exploit strongly convex constraints to obtain $\mathcal O(\varepsilon^{-1/2})$ complexity. Muehlebach and Jordan \cite{MuehlebachMP} accelerate through velocity constraints and local approximations of the feasible set. Without strong convexity, Zhang and Lan \cite{ZhangLanACGD} develop accelerated constrained gradient descent for the smooth-constraint setting and obtain $\mathcal O(k^{-2})$ residual estimates for weighted primal averages under Lipschitz-gradient assumptions. Deng et al. \cite{DengInformsJoc} develop uniformly optimal methods for H\"older-smooth problems.

 The proposed methods give accelerated residual estimates at the generated primal iterates and, in the strict parameter regime, combine $o(k^{-2})$ residual decay with convergence of the  primal-dual sequence.  Table~\ref{tab:comparison} lists representative guarantees under the assumptions of the exists results. Composite means a smooth convex term plus a possibly nonsmooth convex term. Smooth means that the gradients are Lipschitz continuous, $C^1$ means continuously differentiable, and SC means strongly convex.

\subsection{Main contributions}

We develop a class of augmented Lagrangian methods with
Nesterov acceleration for problem~\eqref{ques_main}.
The primal and dual variables use the Nesterov extrapolations
\[
\begin{aligned}
 \bar x_k
 &=x_k+\frac{k-1}{k+\alpha-1}(x_k-x_{k-1}),\\
 \bar\lambda_k
 &=\lambda_k+\frac{k-1}{k+\alpha-1}
   (\lambda_k-\lambda_{k-1}),
 \qquad \alpha\geq3.
\end{aligned}
\]
 The accelerated augmented Lagrangian methods (AALM) combine these extrapolations with difference extrapolation of constraint values and a projected multiplier. Their two convergence regimes are summarized in Table~\ref{tab:comparison}. The main contributions are as follows.
\begin{enumerate}
\item
For $\alpha\geq3$, we establish boundedness of the primal-dual iterates and the simultaneous last-iterate estimates
\[
 |f(x_k)-f^\star|=\mathcal O(k^{-2}),
 \qquad
 \norm{\pos{g(x_k)}}=\mathcal O(k^{-2}).
\]
These bounds include the critical case $\alpha=3$ and allow weighted-summable primal stationarity errors. They require neither the objective strong convexity used for acceleration in \cite{HamedaniAybat2021,ConEx2023,ZhuSiopt} nor the constraint strong convexity used in \cite{LinNeurips}. Our augmented Lagrangian methods use Nesterov acceleration to obtain both residual estimates at the generated primal iterates for composite objectives and convex $C^1$ constraint functions. 

\item
 For $\alpha>3$, we prove convergence of the  primal-dual sequence to a KKT pair and sharpen both estimates to
\[
 |f(x_k)-f^\star|=o(k^{-2}),
 \qquad
 \norm{\pos{g(x_k)}}=o(k^{-2}).
\]
In the general convex regime, ACGD gives $\mathcal O(k^{-2})$ rates for weighted primal averages \cite{ZhangLanACGD}, and smooth APMM gives an $\mathcal O(k^{-2})$ primal merit bound \cite{DengInformsJoc}. The proposed methods establish little-$o$ rates at the generated primal iterates together with full primal-dual convergence. They extend little-$o$ acceleration and iterate convergence from unconstrained \cite{AttouchPtSiopt,ChambolleJota} and affine-constrained \cite{HeArxiv} settings to nonlinear convex inequalities.
\end{enumerate}

\begin{table}[tbp]
\caption{Comparison of convergence results for algorithms solving problem~\eqref{ques_main}}
\label{tab:comparison}
\centering
\footnotesize
\setlength{\tabcolsep}{3pt}
\renewcommand{\arraystretch}{1.12}
\begin{tabularx}{\linewidth}{@{}P{.23\linewidth}YYC{.14\linewidth}C{.17\linewidth}@{}}
\toprule
Method & Objective & Constraints & Rate & \shortstack{Primal--dual\\convergence}\\
\midrule
LALM \cite{XuLinear2021} & Composite & Smooth & $\mathcal O(k^{-1})$ & Yes\\
\midrule
iALM \cite{XuMathProg} & Composite & Smooth & $\mathcal O(k^{-1})$ & --\\
\midrule
\multirow{2}{=}{APD \cite{HamedaniAybat2021}}
& Nonsmooth & Smooth & $\mathcal O(k^{-1})$ & Yes\\
\cmidrule(lr){2-5}
& SC, nonsmooth & Smooth & $\mathcal O(k^{-2})$ & --\\
\midrule
\multirow{2}{=}{ConEx \cite{ConEx2023}}
& Composite & Smooth & $\mathcal O(k^{-1})$ & --\\
\cmidrule(lr){2-5}
& SC, composite & Smooth & $\mathcal O(k^{-2})$ & --\\
\midrule
\multirow{2}{=}{Alg.~1 \cite{ZhuSiopt}}
& Composite & Smooth & $\mathcal O(k^{-1})$ & --\\
\cmidrule(lr){2-5}
& SC, composite & Smooth & $\mathcal O(k^{-2})$ & --\\
\midrule
APDPro \cite{LinNeurips} & Nonsmooth & SC, smooth & $\mathcal O(k^{-2})$ & --\\
\midrule
ACGD~\cite{ZhangLanACGD} & Composite & Smooth & $\mathcal O(k^{-2})$ & --\\
\midrule
APMM \cite{DengInformsJoc} & Smooth & Smooth & $\mathcal O(k^{-2})$ & --\\
\midrule
 AALM ($\alpha\geq3$) & Composite & $C^1$ & $\mathcal O(k^{-2})$ & --\\
\midrule
AALM ($\alpha>3$) & Composite & $C^1$ & $o(k^{-2})$ & Yes\\
\bottomrule
\end{tabularx}
\end{table}

\subsection{Organization}
Section~\ref{sec:algorithm} motivates AALM through primal-dual dynamics, derives its updates, and establishes the energy estimate. Section~\ref{sec:convergence} proves accelerated residual rates and primal-dual convergence. Section~\ref{sec:numerics} reports  numerical experiments. Section~\ref{sec:conclusion} concludes. Appendix~\ref{sec:appendix} collects technical lemmas used in the analysis.

\section{Accelerated augmented Lagrangian methods and energy estimates}\label{sec:algorithm}

Continuous-time interpretations of Nesterov acceleration \cite{SuJmlr,AttouchMathProg} motivate inertial discretizations. For problem~\eqref{ques_main}, we start from the following perturbed primal-dual differential inclusion:
\begin{equation}\label{eq:dy_per}
 \begin{aligned}
 \epsilon(t)&\in\ddot x+\frac{\alpha}{t}\dot x+\nabla\phi(x)+\partial h(x)+J_g(x)^\top p,\\
 0&=\ddot\lambda+\frac{\alpha}{t}\dot\lambda-\frac{\sigma}{\beta}(p-\widehat\lambda),\qquad
 p=\pos{\widehat\lambda+\beta\widehat g},\\
 \widehat\lambda&=\lambda+\frac{t}{\gamma}\dot\lambda,\qquad
 \widehat g=g(x)+\frac{t}{\gamma}J_g(x)\dot x.
 \end{aligned}
\end{equation}
This is a perturbed composite extension of the smooth continuous-time system in \cite{HeNonlinear2026}. We use it to motivate the updates and then analyze the resulting algorithm directly through a discrete energy estimate.

\subsection{Derivation by time discretization}

 The derivation has two steps. We first discretize the inertial terms and the primal-dual coupling. We then eliminate the implicit multiplier, reducing each iteration to one strongly convex primal subproblem and explicit multiplier updates.

Let the fixed time step be $\sqrt{\tau}$ and use the common grid
$
  t_k=(k+\alpha-1)\sqrt{\tau}$, $
  x_k=x(t_k)$, $\lambda_k=\lambda(t_k)$.
Thus $t_0=(\alpha-1)\sqrt{\tau}>0$. At $t_k$, we have
\begin{align*}
  \ddot x(t_k)+\frac{\alpha}{t_k}\dot x(t_k)&\approx \frac{x_{k+1}-2x_k+x_{k-1}}{\tau}
   +\frac{\alpha}{(k+\alpha-1)\sqrt{\tau}}
    \frac{x_k-x_{k-1}}{\sqrt{\tau}}
   =\frac{x_{k+1}-\bar x_k}{\tau},\\
  \ddot \lambda(t_k)+\frac{\alpha}{t_k}\dot \lambda(t_k)&\approx\frac{\lambda_{k+1}-2\lambda_k+\lambda_{k-1}}{\tau}
   +\frac{\alpha}{(k+\alpha-1)\sqrt{\tau}}
    \frac{\lambda_k-\lambda_{k-1}}{\sqrt{\tau}}
   =\frac{\lambda_{k+1}-\bar\lambda_k}{\tau},
 \end{align*}where
\begin{equation}\label{eq:iner}
\begin{aligned}
  (\bar x_k,\bar\lambda_k)
   =( x_k,\lambda_k)+\frac{k-1}{k+\alpha-1}(( x_k,\lambda_k)-(x_{k-1},\lambda_{k-1}))
\end{aligned}
\end{equation}
are the Nesterov extrapolations. For the velocity corrections, we use forward differences over $[t_k,t_{k+1}]$. In particular, the chain rule $\frac{d}{dt}g(x(t))=J_g(x(t))\dot x(t)$ motivates replacing $J_g(x)\dot x$ by $(g(x_{k+1})-g(x_k))/\sqrt\tau$. Thus the constraint correction extrapolates the values of $g$, as specified in part~(d) below. We evaluate $\nabla\phi$ explicitly at $\bar x_k$ and treat $\partial h$ and the multiplier coupling implicitly at $x_{k+1}$. This gives the following  explicit-implicit scheme for \eqref{eq:dy_per}:
\begin{equation}\label{eq:dis}
\left\{
\begin{aligned}
\mathrm{(a)}\quad&
\epsilon_{k+1}\in
\frac{x_{k+1}-\bar x_k}{\tau}
+\nabla\phi(\bar x_k)
+\partial h(x_{k+1})
+J_g(x_{k+1})^\top p_{k+1},
\\[0.25em]
\mathrm{(b)}\quad&
\frac{\lambda_{k+1}-\bar\lambda_k}{\tau}
=
\frac{\sigma}{\beta}
\bigl(p_{k+1}-\widehat\lambda_{k+1}\bigr),
\\[0.25em]
\mathrm{(c)}\quad&
\widehat\lambda_{k+1}
=
\lambda_k+\frac{k+\alpha-1}{\gamma}
(\lambda_{k+1}-\lambda_k),
\\[0.25em]
\mathrm{(d)}\quad&
\widehat g_{k+1}
=
g(x_k)+\frac{k+\alpha-1}{\gamma}
\bigl(g(x_{k+1})-g(x_k)\bigr),
\\[0.25em]
\mathrm{(e)}\quad&
p_{k+1}
=
\pos{
  \widehat\lambda_{k+1}
  +\beta\widehat g_{k+1}
}.
\end{aligned}
\right.
\end{equation}
 Here $\epsilon_{k+1}$ represents the primal stationarity error. Although the primal and multiplier updates in \eqref{eq:dis} are coupled, the projection identity allows us to eliminate $\widehat\lambda_{k+1}$ and express $p_{k+1}$ as a function of the candidate primal point.

To eliminate the implicit corrected multiplier
$\widehat\lambda_{k+1}$, define
\begin{equation}\label{eq:dual}
  \widetilde\lambda_k
  :=
  \lambda_k+\frac{k-1}{\gamma}
  (\lambda_k-\lambda_{k-1}),
  \qquad
  c_{k+1}
  :=
  \beta+\frac{\sigma\tau(k+\alpha-1)}{\gamma}.
\end{equation}
Substituting \eqref{eq:dis}\textnormal{(c)} into
\eqref{eq:dis}\textnormal{(b)} gives
\[
  \lambda_{k+1}
  =
  \lambda_k
  +\frac{\beta}{c_{k+1}}
  (\bar\lambda_k-\lambda_k)
  +\frac{\sigma\tau}{c_{k+1}}
  (p_{k+1}-\lambda_k).
\]
By \eqref{eq:iner},
\eqref{eq:dis}\textnormal{(c)}, and
\eqref{eq:dual},
\[
\begin{aligned}
  \widehat\lambda_{k+1}
  &=
  \lambda_k
  +
  \frac{\beta}{c_{k+1}}
  \frac{k+\alpha-1}{\gamma}
  (\bar\lambda_k-\lambda_k)
  +
  \frac{\sigma\tau(k+\alpha-1)}
     {\gamma c_{k+1}}
  (p_{k+1}-\lambda_k)\\
  &=
  \frac{\beta}{c_{k+1}}\widetilde\lambda_k
  +
  \left(
    1-\frac{\beta}{c_{k+1}}
  \right)p_{k+1}.
\end{aligned}
\]
Substituting this identity into
\eqref{eq:dis}\textnormal{(e)} gives
\[
	 p_{k+1}
  =
  \proj_{\R_+^m}
  \left(
    \frac{\beta}{c_{k+1}}
    \bigl(
      \widetilde\lambda_k
      +c_{k+1}\widehat g_{k+1}
    \bigr)
    +
    \left(
      1-\frac{\beta}{c_{k+1}}
    \right)p_{k+1}
  \right).
\]
Since $\beta/c_{k+1}\in(0,1)$, Lemma~\ref{lem:proj} with
$\theta=\beta/c_{k+1}$ yields
\begin{equation}\label{eq:def_pk}
  p_{k+1}
  =
  \proj_{\R_+^m}
  \bigl(
    \widetilde\lambda_k
    +c_{k+1}\widehat g_{k+1}
  \bigr)
  =
  \pos{
    \widetilde\lambda_k
    +c_{k+1}\widehat g_{k+1}
  }.
\end{equation}

The right-hand side of \eqref{eq:def_pk} depends on $x_{k+1}$ only through the extrapolated constraint values. To express the primal update as an optimization subproblem, define, for $x\in\R^n$,
\begin{equation}\label{eq:dis_ghat}
  \widehat g_{k+1}(x)
  :=
  g(x_k)
  +\frac{k+\alpha-1}{\gamma}
  \bigl(g(x)-g(x_k)\bigr)
\end{equation}
and
\begin{equation}\label{eq:Theta}
\begin{aligned}
  \Theta_{k+1}(x)
  :={}&
  h(x)
  +\frac1{2\tau}
  \norm{
    x-\bigl(\bar x_k-\tau\nabla\phi(\bar x_k)\bigr)
  }^2\\
  &+
  \frac{\gamma}
     {2c_{k+1}(k+\alpha-1)}
  \norm{
    \pos{
      \widetilde\lambda_k
      +c_{k+1}\widehat g_{k+1}(x)
    }
  }^2.
\end{aligned}
\end{equation}
Each component of $\widehat g_{k+1}$ is convex, and $s\mapsto[s]_+^2$ is convex and nondecreasing, so the penalty in \eqref{eq:Theta} is convex. The quadratic term makes $\Theta_{k+1}$ $1/\tau$-strongly convex. Since $h$ has an affine lower bound and the penalty is nonnegative, $\Theta_{k+1}$ has a unique minimizer. Differentiating its smooth part gives
\begin{equation}\label{eq:part_theta}
    \partial\Theta_{k+1}(x)
    =
    \partial h(x)
    +\frac1\tau
    \bigl[x-(\bar x_k-\tau\nabla\phi(\bar x_k))\bigr]
    +J_g(x)^\top
    \pos{\widetilde\lambda_k+c_{k+1}\widehat g_{k+1}(x)}.
\end{equation}
 Motivated by the perturbation in
\eqref{eq:dis}\textnormal{(a)}, we solve the
primal subproblem inexactly. Let
$\{\varepsilon_k\}_{k\geq2}\subset\R_+$. At iteration $k$, we compute
$x_{k+1}\in\operatorname{dom}h$ satisfying
\begin{equation}\label{eq:sub}
  \operatorname{dist}
  \bigl(
    0,\partial\Theta_{k+1}(x_{k+1})
  \bigr)
  \leq\varepsilon_{k+1}.
\end{equation}
By \eqref{eq:sub}, there exists
$\epsilon_{k+1}\in\partial\Theta_{k+1}(x_{k+1})$ such that
$\norm{\epsilon_{k+1}}\leq\varepsilon_{k+1}$.
Using \eqref{eq:def_pk} and \eqref{eq:part_theta}, this choice of
$\epsilon_{k+1}$ yields exactly the perturbed primal inclusion
\eqref{eq:dis}\textnormal{(a)}. Let $\bar{x}=\arg\min_{x\in\R^n}\Theta_{k+1}(x)$ denote the unique exact minimizer of
$\Theta_{k+1}(x)$. Since $\Theta_{k+1}(x)$ is $1/\tau$-strongly convex and $\epsilon_{k+1}\in\partial\Theta_{k+1}(x_{k+1})$,
\[
  \Theta_{k+1}(\bar{x})
  \geq
  \Theta_{k+1}(x_{k+1})
  +
  \ip{\epsilon_{k+1}}
    {\bar{x}-x_{k+1}}
  +
  \frac{1}{2\tau}
  \norm{\bar{x}-x_{k+1}}^2.
\]
Hence, using $\norm{\epsilon_{k+1}}\leq\varepsilon_{k+1}$ and
$ab\leq\frac{1}{2\tau}a^2+\frac{\tau}{2}b^2$
  for all $a,b\geq0$,
we obtain
\[
\begin{aligned}
  \Theta_{k+1}(x_{k+1})
  -\Theta_{k+1}(\bar{x})
  &\leq
  \norm{\epsilon_{k+1}}
  \norm{x_{k+1}-\bar{x}}
  -
  \frac{1}{2\tau}
  \norm{x_{k+1}-\bar{x}}^2\leq
  \frac{\tau}{2}\varepsilon_{k+1}^2.
\end{aligned}
\]
Therefore,
\[
  \Theta_{k+1}(x_{k+1})
  \leq
  \min_{x\in\R^n}\Theta_{k+1}(x)
  +
  \frac{\tau}{2}\varepsilon_{k+1}^2,
\]
so the stationarity test \eqref{eq:sub} also gives a $\frac{\tau}{2}\varepsilon_{k+1}^2$ objective-gap guarantee. Setting $\varepsilon_{k+1}=0$ recovers the exact minimizer.

These identities give Algorithm~\ref{alg:inexact}: form the extrapolations, solve the primal subproblem, evaluate the projected multiplier, and update the dual variable. Algorithm~\ref{alg:inexact} and \eqref{eq:dis} are equivalent in the following sense. Every primal update satisfying \eqref{eq:sub} admits an error $\epsilon_{k+1}\in\partial\Theta_{k+1}(x_{k+1})$ with $\norm{\epsilon_{k+1}}\leq\varepsilon_{k+1}$; the explicit multiplier updates then satisfy \eqref{eq:dis}\textnormal{(b)--(e)}. Conversely, eliminating the corrected multiplier from \eqref{eq:dis} gives the same primal residual test and multiplier updates whenever $\norm{\epsilon_{k+1}}\leq\varepsilon_{k+1}$. We therefore use \eqref{eq:dis} to analyze the iterates of Algorithm~\ref{alg:inexact}.

\begin{algorithm}[htbp]
\caption{Accelerated augmented Lagrangian method (AALM)}
\label{alg:inexact}
\begin{algorithmic}[1]
\Require
$x_0=x_1\in\operatorname{dom}h$,
$\lambda_0=\lambda_1\in\R^m$;
$\alpha\geq3$, $2\leq\gamma\leq\alpha-1$,
$0<\tau\leq1/L_\phi$, $\sigma>0$, $\beta>0$;
$\{\varepsilon_k\}_{k\geq2}\subset\R_+$.
\For{$k=1,2,\ldots$}
  \State Compute
  $\displaystyle
    \bar x_k\gets
    x_k+\frac{k-1}{k+\alpha-1}(x_k-x_{k-1})$.

  \State Compute
  $\displaystyle
    \bar\lambda_k\gets
    \lambda_k+\frac{k-1}{k+\alpha-1}
    (\lambda_k-\lambda_{k-1})$
  and
  $\displaystyle
    \widetilde\lambda_k\gets
    \lambda_k+\frac{k-1}{\gamma}
    (\lambda_k-\lambda_{k-1})$.

  \State Set
  $\displaystyle
    c_{k+1}\gets
    \beta+\frac{\sigma\tau(k+\alpha-1)}{\gamma}$, and define $\widehat g_{k+1}(x)$ and $\Theta_{k+1}(x)$ by \eqref{eq:dis_ghat} and \eqref{eq:Theta}, respectively.

  \State Find $x_{k+1}\in\operatorname{dom}h$ such that
  \eqref{eq:sub} holds. Then
  $x_{k+1}$ is a
  $\frac{\tau}{2}\varepsilon_{k+1}^2$-optimal solution of
  $\min_{x\in\R^n}\Theta_{k+1}(x)$, i.e.
  \[
    \Theta_{k+1}(x_{k+1})
    \leq
    \min_{x\in\R^n}\Theta_{k+1}(x)
    +\frac{\tau}{2}\varepsilon_{k+1}^2.
  \]

  \State Set
  $\displaystyle
    \widehat g_{k+1}\gets
    \widehat g_{k+1}(x_{k+1})$
  and
  $\displaystyle
    p_{k+1}\gets
    \pos{
      \widetilde\lambda_k
      +c_{k+1}\widehat g_{k+1}
    }$.

  \State Update
  $\displaystyle
    \lambda_{k+1}\gets
    \lambda_k
    +\frac{\beta}{c_{k+1}}
    (\bar\lambda_k-\lambda_k)
    +\frac{\sigma\tau}{c_{k+1}}
    (p_{k+1}-\lambda_k)$.
\EndFor
\end{algorithmic}
\end{algorithm}

\begin{remark}
To identify AALM's augmented Lagrangian structure, recall the Powell-Hestenes-Rockafellar (PHR) function, for $\rho>0$ and $\nu\in\R^m$:
\[ \mathscr A_\rho(x;\nu)=f(x)+\frac{\norm{\pos{\nu+\rho g(x)}}^2-\norm\nu^2}{2\rho};
\]
see \cite{Powell1969,Hestenes1969,RockALM1976}. Set $a_k=(k+\alpha-1)/\gamma$, $\rho_{k+1}=c_{k+1}a_k$, and
\[ \nu_k=\widetilde\lambda_k+c_{k+1}(1-a_k)g(x_k).\]
Then $\widetilde\lambda_k+c_{k+1}\widehat g_{k+1}(x)=\nu_k+\rho_{k+1}g(x)$, and~\eqref{eq:Theta} becomes
\begin{equation*}
\begin{aligned}
 \Theta_{k+1}(x)={}&h(x)+\ip{\nabla\phi(\bar x_k)}{x-\bar x_k}
 +\frac{\norm{x-\bar x_k}^2}{2\tau}\\
 &+\frac{\norm{\pos{\nu_k+\rho_{k+1}g(x)}}^2-\norm{\nu_k}^2}{2\rho_{k+1}}
 +\text{constant},
\end{aligned}
\end{equation*}
The constant is independent of $x$. Thus AALM uses a proximal, partially linearized PHR subproblem with the shifted multiplier $\nu_k$. Relative to classical multiplier methods \cite{Powell1969,Hestenes1969,RockALM1976}, the Nesterov extrapolations and the constraint-dependent shift provide the inertial coupling used in the accelerated energy estimate; $p_{k+1}$ enters the primal subproblem, while $\lambda_{k+1}$ follows the explicit update in Algorithm~\ref{alg:inexact}.
\end{remark}

\begin{remark}
The primal subproblem is $1/\tau$-strongly convex.
Suppose that the proximal map of $h$ is computable. If $g$ is smooth, the differentiable part of $\Theta_{k+1}$ has a Lipschitz continuous gradient, then we can use standard FISTA \cite{BeckSiis} or a variant that exploits strong convexity \cite{Beck2017} with a fixed stepsize. If $J_g$ is only locally Lipschitz continuous, the adaptive accelerated proximal-gradient method in \cite[Algorithm~2]{LuMeiSiopt} can choose the inner stepsizes by backtracking. Universal methods are an option under H\"older continuity of the differentiable part's gradient \cite{Nesterov15}. The inner iterations terminate when \eqref{eq:sub} is satisfied.
The outer convergence analysis requires only continuous differentiability of $g$, relaxing the Lipschitz smoothness requirements on the constraints in \cite{ZhangLanACGD,HamedaniAybat2021,ZhuSiopt,LinNeurips}.
\end{remark}

\subsection{Energy estimate}

We now construct the discrete energy sequence and establish its one-step estimate.

 Let $\{(x_k,\lambda_k)\}_{k\geq0}$ be generated by
Algorithm~\ref{alg:inexact}, and let
$(x^\star,\lambda^\star)\in\Omega$ be a KKT pair of
\eqref{ques_main}. Set $z^\star=(x^\star,\lambda^\star)$. For $k\geq1$, define
\begin{equation}\label{eq:aux}
\begin{aligned}
  \ell_k
  &:=
  \mathcal L(x_k,\lambda^\star)
  -\mathcal L(x^\star,\lambda^\star),\\
  U_k
  &:=
  \gamma(x_k-x^\star)
  +(k-1)(x_k-x_{k-1}),\\
  V_k
  &:=
  \gamma(\lambda_k-\lambda^\star)
  +(k-1)(\lambda_k-\lambda_{k-1}).
\end{aligned}
\end{equation}
The discrete energy sequence is defined by
\begin{equation}\label{eq:dis_energy}
\begin{aligned}
  {\mathcal E_{k,z^\star}}:={}&
  \tau r_{k-1}^2\ell_k
  +\frac12\norm{U_k}^2
  +\frac{\gamma\delta}{2}\norm{x_k-x^\star}^2+\frac1{2\sigma}\norm{V_k}^2
  +\frac{\gamma\delta}{2\sigma}
   \norm{\lambda_k-\lambda^\star}^2,
\end{aligned}
\end{equation}
where
\begin{equation}\label{eq_rk}
	 r_k = k+\alpha-1,\quad \text{and}\quad \delta=\alpha-\gamma-1\geq0.
\end{equation}
By the saddle-point property \eqref{eq:saddle}, $\ell_k\geq0$.
Consequently, ${\mathcal E_{k,z^\star}}\geq0$ for every $k\geq1$.

The following lemma provides the one-step estimate for
the energy sequence.

\begin{lemma}
\label{lem:dis_est}
Let $\{(x_k,\lambda_k)\}_{k\geq0}$ be generated by
Algorithm~\ref{alg:inexact}, and let
$(x^\star,\lambda^\star)\in\Omega$. For every $k\geq1$,
\begin{align}
 {\mathcal E_{k+1,z^\star}-\mathcal E_{k,z^\star}}
 \leq{}&
 -\tau\bigl((\gamma-2)r_k+1\bigr)\ell_k
 \notag\\
 &
 -\frac{\delta}{2}(2k+\gamma+\delta)
 \left(
   \|x_{k+1}-x_k\|^2
   +\frac1\sigma\|\lambda_{k+1}-\lambda_k\|^2
 \right)
 \notag\\
 &-\frac{\tau\gamma r_k}{\beta}
 \|p_{k+1}-\widehat\lambda_{k+1}\|^2
 +\tau r_k
 \ip{\epsilon_{k+1}}
   {U_{k+1}+\delta(x_{k+1}-x_k)}.
 \label{eq:dis_full_est}
\end{align}
\end{lemma}

\begin{proof}
By the inexact primal update
\eqref{eq:dis}\textnormal{(a)}, there exists
$\xi_{k+1}\in\partial h(x_{k+1})$ such that
\begin{equation}\label{eq:res}
  \epsilon_{k+1}
  =
  \frac{x_{k+1}-\bar x_k}{\tau}
  +\nabla\phi(\bar x_k)
  +\xi_{k+1}
  +J_g(x_{k+1})^\top p_{k+1}.
\end{equation}
Let $x\in\dom h$. The $L_\phi$-Lipschitz continuity of
$\nabla\phi$ gives
\[
  \phi(x_{k+1})
  \leq
  \phi(\bar x_k)
  +\ip{\nabla\phi(\bar x_k)}{x_{k+1}-\bar x_k}
  +\frac{L_\phi}{2}
   \norm{x_{k+1}-\bar x_k}^2,
\]
whereas the convexity of $\phi$ gives
$
  \phi(x)
  \geq
  \phi(\bar x_k)
  +\ip{\nabla\phi(\bar x_k)}{x-\bar x_k}.
$
Therefore,
\begin{equation}\label{eq:phi_est}
  \phi(x_{k+1})-\phi(x)
  \leq
  \ip{\nabla\phi(\bar x_k)}{x_{k+1}-x}
  +\frac{L_\phi}{2}
   \norm{x_{k+1}-\bar x_k}^2.
\end{equation}
Since $\xi_{k+1}\in\partial h(x_{k+1})$, we also have
\begin{equation}\label{eq:h_est}
  h(x_{k+1})-h(x)
  \leq
  \ip{\xi_{k+1}}{x_{k+1}-x}.
\end{equation}
Moreover, the convexity of $g_i$ gives
$
  g_i(x)
  \geq
  g_i(x_{k+1})
  +\ip{\nabla g_i(x_{k+1})}{x-x_{k+1}}.
$
Multiplying by $p_{k+1,i}\geq0$ and summing over
$i=1,\ldots,m$ yields
\begin{equation}\label{eq:g_est}
  \ip{p_{k+1}}{g(x_{k+1})-g(x)}
  \leq
  \ip{J_g(x_{k+1})^\top p_{k+1}}
    {x_{k+1}-x}.
\end{equation}

Adding \eqref{eq:phi_est}, \eqref{eq:h_est}, and \eqref{eq:g_est} and
using $f=\phi+h$, we obtain, for every $x\in\dom h$,
\begin{align}\label{eq:f_est}
  &f(x_{k+1})+\ip{p_{k+1}}{g(x_{k+1})} -f(x)-\ip{p_{k+1}}{g(x)}
  \notag\\
  &\quad\leq
  \ip{
    \nabla\phi(\bar x_k)
    +\xi_{k+1}
    +J_g(x_{k+1})^\top p_{k+1}
  }{
    x_{k+1}-x
     }+\frac{L_\phi}{2}
   \norm{x_{k+1}-\bar x_k}^2
  \\
  &\quad\overset{\eqref{eq:res}}{=}
  \ip{\epsilon_{k+1}}{x_{k+1}-x}
  -\frac1\tau
   \ip{x_{k+1}-\bar x_k}{x_{k+1}-x}+\frac{L_\phi}{2}
   \norm{x_{k+1}-\bar x_k}^2
  \notag\\
  &\quad=
  \frac1{2\tau}
  \left(
    \norm{\bar x_k-x}^2
    -\norm{x_{k+1}-x}^2
  \right)
  +\ip{\epsilon_{k+1}}{x_{k+1}-x}
  -\frac{1-\tau L_\phi}{2\tau}
   \norm{x_{k+1}-\bar x_k}^2
  \notag\\
  &\quad\leq
  \frac1{2\tau}
  \left(
    \norm{\bar x_k-x}^2
    -\norm{x_{k+1}-x}^2
  \right)
  +\ip{\epsilon_{k+1}}{x_{k+1}-x},\notag
\end{align}
where the last identity follows from $-2\ip{a}{b}
  =
  \norm{a-b}^2-\norm{a}^2-\norm{b}^2
$. The condition $\tau\leq1/L_\phi$ makes the last inequality valid.
Apply \eqref{eq:f_est} with $x=x^\star$ and $x=x_k$, and multiply the
resulting inequalities by $\tau\gamma r_k$ and
$\tau r_k(r_k-\gamma)$, respectively. Adding the two weighted
inequalities, we first simplify their left-hand side. By the definition of
$\widehat g_{k+1}$ in \eqref{eq:dis}(d),
\[
  r_k g(x_{k+1})
  -(r_k-\gamma)g(x_k)
  =
  \gamma\widehat g_{k+1}.
\]
Using this identity, the definition of $\ell_k$, and
$\ip{\lambda^\star}{g(x^\star)}=0$, the sum of the two weighted
left-hand sides, namely, $\tau\gamma r_k$ times \eqref{eq:f_est} with
$x=x^\star$ plus $\tau r_k(r_k-\gamma)$ times \eqref{eq:f_est} with
$x=x_k$, can be written as
\begin{align*}
  \tau r_k
  \Big[
    &r_k\ell_{k+1}
    -(r_k-\gamma)\ell_k
    +\gamma
    \ip{p_{k+1}-\lambda^\star}{\widehat g_{k+1}}
    -\gamma\ip{p_{k+1}}{g(x^\star)}
  \Big].
\end{align*}
Since $g(x^\star)\leq0$ and $p_{k+1}\geq0$, $
  \ip{p_{k+1}}{g(x^\star)}\leq0.
$
Moreover, since $r_{k-1}=r_k-1$,
\begin{equation}\label{eq_rk_eqa}
  r_k(r_k-\gamma)
  =
  r_{k-1}^2-\bigl((\gamma-2)r_k+1\bigr).
\end{equation}

It follows that the sum of the two weighted left-hand sides is bounded
below by
\begin{align*}
  &\tau r_k^2\ell_{k+1}
  -\tau r_{k-1}^2\ell_k
  +\tau\bigl((\gamma-2)r_k+1\bigr)\ell_k
  +\tau\gamma r_k
    \ip{p_{k+1}-\lambda^\star}{\widehat g_{k+1}}.
\end{align*}
Therefore, adding the two weighted inequalities obtained from
\eqref{eq:f_est} gives
\begin{align}\label{eq_left_est}
  &\tau r_k^2\ell_{k+1}
  -\tau r_{k-1}^2\ell_k
  +\tau\bigl(
    (\gamma-2)r_k+1
  \bigr)\ell_k
     +\tau\gamma r_k
    \ip{p_{k+1}-\lambda^\star}{\widehat g_{k+1}}
  \notag\\
  &\qquad \leq
  \frac{r_k}{2}
  \Big[
    \gamma
    \bigl(
      \norm{\bar x_k-x^\star}^2
      -\norm{x_{k+1}-x^\star}^2
    \bigr)
  +(r_k-\gamma)
    \bigl(
      \norm{\bar x_k-x_k}^2
      -\norm{x_{k+1}-x_k}^2
    \bigr)
  \Big]
  \notag\\
  &\qquad\quad+
  \tau r_k
  \ip{\epsilon_{k+1}}
    {\gamma(x_{k+1}-x^\star)
    +(r_k-\gamma)(x_{k+1}-x_k)}.
\end{align}

We next simplify the quadratic terms on the right-hand side of
\eqref{eq_left_est}. By the definition of $\bar x_k$,
$
  \bar x_k-x^\star
  =
  (x_k-x^\star)
  +\frac{k-1}{r_k}(x_k-x_{k-1}),
$
$
  \bar x_k-x_k
  =
  \frac{k-1}{r_k}(x_k-x_{k-1}),
$
and
$
  x_{k+1}-x^\star
  =
  (x_k-x^\star)+(x_{k+1}-x_k).
$
Hence,
\begin{align*}
  &\frac12
  \left(
    \norm{\bar x_k-x^\star}^2
    -\norm{x_{k+1}-x^\star}^2
  \right)
  \nonumber\\
  &\quad=
  \frac12
  \left(
    \norm{x_k-x^\star+\tfrac{k-1}{r_k}(x_k-x_{k-1})}^2
    -\norm{(x_k-x^\star)+(x_{k+1}-x_k)}^2
  \right)
  \notag\\
  &\quad=
  \frac{k-1}{r_k}
  \ip{x_k-x^\star}{x_k-x_{k-1}}
  -
  \ip{x_k-x^\star}{x_{k+1}-x_k}
  \notag\\
  &\qquad+
  \frac{(k-1)^2}{2r_k^2}
  \norm{x_k-x_{k-1}}^2
  -
  \frac12\norm{x_{k+1}-x_k}^2,
  \end{align*}
and
\[
  \frac12
  \left(
    \norm{\bar x_k-x_k}^2
    -\norm{x_{k+1}-x_k}^2
  \right)
  =
  \frac{(k-1)^2}{2r_k^2}
  \norm{x_k-x_{k-1}}^2
  -
  \frac12\norm{x_{k+1}-x_k}^2.
\]
 On the other hand, by the definition of $U_k$ and \eqref{eq_rk},
\[
  \frac12\norm{U_k}^2
  =
  \frac{\gamma^2}{2}\norm{x_k-x^\star}^2
  +\gamma(k-1)
  \ip{x_k-x^\star}{x_k-x_{k-1}}+\frac{(k-1)^2}{2}
  \norm{x_k-x_{k-1}}^2,
\]
and\[
\begin{aligned}
  &\frac12
  \norm{
    U_{k+1}+\delta(x_{k+1}-x_k)
  }^2=\frac12
  \norm{
     \gamma(x_k-x^\star)+r_k(x_{k+1}-x_k)
  }^2\\
  &=
  \frac{\gamma^2}{2}\norm{x_k-x^\star}^2
  +\gamma r_k
  \ip{x_k-x^\star}{x_{k+1}-x_k}
  +\frac{r_k^2}{2}
  \norm{x_{k+1}-x_k}^2.
\end{aligned}
\]
Combining these four identities yields
\begin{align}\label{eq:es_2}
  &\frac{r_k}{2}
  \Big[
    \gamma
    \bigl(
      \norm{\bar x_k-x^\star}^2
      -\norm{x_{k+1}-x^\star}^2
    \bigr)
    +(r_k-\gamma)
    \bigl(
      \norm{\bar x_k-x_k}^2
      -\norm{x_{k+1}-x_k}^2
    \bigr)
  \Big] \notag
  \\
  &\quad=
  \gamma(k-1)
  \ip{x_k-x^\star}{x_k-x_{k-1}}
  -
  \gamma r_k
  \ip{x_k-x^\star}{x_{k+1}-x_k}
 \notag
 \\
  &\qquad+
  \frac{(k-1)^2}{2}
  \norm{x_k-x_{k-1}}^2
  -
  \frac{r_k^2}{2}
  \norm{x_{k+1}-x_k}^2 \\
    &\quad=
  \frac12\norm{U_k}^2
  -
  \frac12
  \norm{
    U_{k+1}+\delta(x_{k+1}-x_k)
  }^2.\notag
\end{align}

It remains to express the last term in
\eqref{eq:es_2} in terms of the primal
energy. From \eqref{eq:aux}, we have
\begin{align*}
  &\frac12
  \norm{
    U_{k+1}+\delta(x_{k+1}-x_k)
  }^2=\frac12\norm{U_{k+1}}^2
  +
  \delta\ip{U_{k+1}}{x_{k+1}-x_k}
  +
  \frac{\delta^2}{2}
  \norm{x_{k+1}-x_k}^2
  \notag\\
  &=
  \frac12\norm{U_{k+1}}^2
  +
  \gamma\delta
  \ip{x_{k+1}-x^\star}{x_{k+1}-x_k}
 +\left(
    k\delta+\frac{\delta^2}{2}
  \right)
  \norm{x_{k+1}-x_k}^2\\
  &=
  \frac12\norm{U_{k+1}}^2
  +
  \frac{\gamma\delta}{2}
  \left(
    \norm{x_{k+1}-x^\star}^2
    -\norm{x_k-x^\star}^2
  \right)+\frac{\delta}{2}
  (2k+\gamma+\delta)
  \norm{x_{k+1}-x_k}^2.\nonumber
 \end{align*}
 Substituting this into
\eqref{eq:es_2} and using \eqref{eq_left_est}, we obtain
\begin{align} \label{eq:dis_pri}
  &\tau r_k^2\ell_{k+1}
  -\tau r_{k-1}^2\ell_k
  +\frac12
  \left(
    \norm{U_{k+1}}^2-\norm{U_k}^2
  \right)\notag\\
  &\quad+
  \frac{\gamma\delta}{2}
  \left(
    \norm{x_{k+1}-x^\star}^2
    -\norm{x_k-x^\star}^2
  \right)+
  \frac{\delta}{2}(2k+\gamma+\delta)
  \norm{x_{k+1}-x_k}^2
 \notag\\
  &\ \leq
  -\tau\bigl((\gamma-2)r_k+1\bigr)\ell_k
  -
  \tau\gamma r_k
  \ip{p_{k+1}-\lambda^\star}{\widehat g_{k+1}}\notag\\
  &\qquad+
  \tau r_k
  \ip{\epsilon_{k+1}}
    {U_{k+1}+\delta(x_{k+1}-x_k)}.
\end{align}

For the dual part, by the definitions of $\bar\lambda_k$,
$\widehat\lambda_{k+1}$, and $V_k$, together with
\eqref{eq:dis}\textnormal{(b)} and \eqref{eq:dis}\textnormal{(c)}, we have
\begin{align*}
  V_{k+1}-V_k+\delta(\lambda_{k+1}-\lambda_k)
  &=
  r_k(\lambda_{k+1}-\bar\lambda_k)=
  \frac{\tau\sigma r_k}{\beta}
  \bigl(p_{k+1}-\widehat\lambda_{k+1}\bigr),
  \\
  V_{k+1}+\delta(\lambda_{k+1}-\lambda_k)
  &=
  \gamma(\widehat\lambda_{k+1}-\lambda^\star).
  \end{align*}
Together with the identity
$
  \ip{a-b}{a}
  =
  \frac12\bigl(
    \norm{a}^2-\norm{b}^2+\norm{a-b}^2
  \bigr),
$
these relations imply
\begin{align}\label{eq:es3}
  &\frac1{2\sigma}
  \left(
    \norm{
      V_{k+1}+\delta(\lambda_{k+1}-\lambda_k)
    }^2
    -\norm{V_k}^2+\norm{
    V_{k+1}-V_k+\delta(\lambda_{k+1}-\lambda_k)
  }^2
  \right)
  \notag\\
 &\qquad=\frac1\sigma
  \ip{
    V_{k+1}-V_k+\delta(\lambda_{k+1}-\lambda_k)
  }{
    V_{k+1}+\delta(\lambda_{k+1}-\lambda_k)
  }
  \\
  &\qquad=
  \frac{\tau\gamma r_k}{\beta}
  \ip{\widehat\lambda_{k+1}-\lambda^\star}
    {p_{k+1}-\widehat\lambda_{k+1}}.\notag
\end{align}
Moreover, a computation analogous to that for $U_k$ gives the exact
identity
\begin{align*}
  \frac1{2\sigma}
  \norm{
    V_{k+1}+\delta(\lambda_{k+1}-\lambda_k)
  }^2
=  &
  \frac1{2\sigma}\norm{V_{k+1}}^2
  +
  \frac{\gamma\delta}{2\sigma}
  \left(
    \norm{\lambda_{k+1}-\lambda^\star}^2
    -\norm{\lambda_k-\lambda^\star}^2
  \right)
  \\
  &+
  \frac{\delta}{2\sigma}(2k+\gamma+\delta)
  \norm{\lambda_{k+1}-\lambda_k}^2.
\end{align*}
Substituting this into
\eqref{eq:es3} yields
\begin{align}\label{eq:dis_dual}
  &\frac1{2\sigma}
  \bigl(
    \norm{V_{k+1}}^2-\norm{V_k}^2
  \bigr)
  +
  \frac{\gamma\delta}{2\sigma}
  \left(
    \norm{\lambda_{k+1}-\lambda^\star}^2
    -\norm{\lambda_k-\lambda^\star}^2
  \right)
+\frac{\delta}{2\sigma}(2k+\gamma+\delta)
  \norm{\lambda_{k+1}-\lambda_k}^2
  \notag\\
   &\qquad\leq
  \frac{\tau\gamma r_k}{\beta}
  \ip{\widehat\lambda_{k+1}-\lambda^\star}
    {p_{k+1}-\widehat\lambda_{k+1}}.
\end{align}

Adding \eqref{eq:dis_pri} and
\eqref{eq:dis_dual}, and using the definition
of ${\mathcal E_{k,z^\star}}$ in \eqref{eq:dis_energy}, we obtain
\begin{align}\label{eq:dis_enerEs}
 {\mathcal E_{k+1,z^\star}-\mathcal E_{k,z^\star}}
 &  \leq -\tau\bigl((\gamma-2)r_k+1\bigr)\ell_k+
  \tau r_k
  \ip{\epsilon_{k+1}}
    {U_{k+1}+\delta(x_{k+1}-x_k)}
  \notag\\
  &-
  \frac{\delta}{2}(2k+\gamma+\delta)
  \left(
    \norm{x_{k+1}-x_k}^2
    +\frac1\sigma\norm{\lambda_{k+1}-\lambda_k}^2
  \right)
  \notag\\
  &
  +\tau\gamma r_k
  \Big[
    -\ip{\widehat g_{k+1}}
      {p_{k+1}-\lambda^\star}
    +
    \frac1\beta
    \ip{
      \widehat\lambda_{k+1}-\lambda^\star
    }{
      p_{k+1}-\widehat\lambda_{k+1}
    }
  \Big].
 \end{align}
It remains to estimate the last two coupling terms. By
\eqref{eq:dis}\textnormal{(e)},
$
  p_{k+1}
  =
  \proj_{\R_+^m}
  \bigl(
    \widehat\lambda_{k+1}
    +\beta\widehat g_{k+1}
  \bigr).
$
Since $\lambda^\star\in\R_+^m$, Lemma~\ref{lem:proj} gives
\[
  \ip{
    \widehat\lambda_{k+1}
    +\beta\widehat g_{k+1}
    -p_{k+1}
  }{
    \lambda^\star-p_{k+1}
  }
  \leq0.
\]
Furthermore,
\begin{align}
  &-\ip{\widehat g_{k+1}}
      {p_{k+1}-\lambda^\star}
  +
  \frac1\beta
  \ip{
    \widehat\lambda_{k+1}-\lambda^\star
  }{
    p_{k+1}-\widehat\lambda_{k+1}
  }
  \notag\\
  &\qquad=
  \frac1\beta
  \ip{
    \widehat\lambda_{k+1}
    +\beta\widehat g_{k+1}
    -p_{k+1}
  }{
    \lambda^\star-p_{k+1}
  }
-\frac1\beta
  \norm{
    p_{k+1}-\widehat\lambda_{k+1}
  }^2
  \notag\\
  &\qquad\leq
  -\frac1\beta
  \norm{
    p_{k+1}-\widehat\lambda_{k+1}
  }^2.\nonumber
\end{align}
Substituting this estimate into \eqref{eq:dis_enerEs} yields the result.
\end{proof}

\section{Convergence analysis}\label{sec:convergence}
In this section, we first use the energy estimate to prove boundedness and simultaneous objective and feasibility rates. Under the strict parameter conditions, the additional dissipation identifies all cluster points as KKT pairs and yields convergence of the full sequence together with little-$o$ refinements.
\subsection{Convergence rates of the nonlinear feasibility and
the objective residual}
The preceding one-step estimate yields the following conclusions under weighted summability of the primal stationarity errors.
\begin{theorem}\label{thm:dis_rates}
Suppose that $\sum_{k=2}^{+\infty}k\varepsilon_k<+\infty$, and let
$\{(x_k,\lambda_k)\}_{k\geq0}$ be generated by
Algorithm~\ref{alg:inexact}. Then the following statements hold:
\begin{itemize}
  \item[(i)] The sequence $\{(x_k,\lambda_k)\}$ is bounded, and
\[
    \|x_k-x_{k-1}\|+\|\lambda_k-\lambda_{k-1}\|
    =\mathcal O\left(\frac{1}{k}\right).
\]
  \item[(ii)] The constraint violation and objective residual satisfy
	\[ \|\pos{g(x_k)}\| =   \mathcal O\left(\frac{1}{k^{2}}\right), \qquad |f(x_k)-f^\star| = \mathcal O\left(\frac{1}{k^{2}}\right). \]
\end{itemize}
\end{theorem}
 \begin{proof}
(i) The definition of ${\mathcal E_{k,z^\star}}$ in \eqref{eq:dis_energy} gives
\[
  \tfrac12\norm{U_{k+1}}^2
  +
  \tfrac{\gamma\delta}{2}
  \norm{x_{k+1}-x^\star}^2
  \leq
  {\mathcal E_{k+1,z^\star}}.
\]
For $k\geq1$, since 
\[
  U_{k+1}+\delta(x_{k+1}-x_k)
  =
  \left(1+\frac{\delta}{k}\right)U_{k+1}
  -\frac{\delta\gamma}{k}(x_{k+1}-x^\star),
\] there exists a constant $C>0$, independent of $k$, such that
\[
  \norm{
    U_{k+1}+\delta(x_{k+1}-x_k)
  }^2
  \leq
  C{\mathcal E_{k+1,z^\star}}.
\]
Consequently, using
$\norm{\epsilon_{k+1}}\leq\varepsilon_{k+1}$ and Lemma~\ref{lem:dis_est}, we obtain
\[\norm{
    U_{k+1}+\delta(x_{k+1}-x_k)
  }^2 \le C{\mathcal E_{1,z^\star}}
  +
  \tau C
  \sum_{j=1}^{k}
  r_j\varepsilon_{j+1}\|U_{j+1}+\delta(x_{j+1}-x_j)\|.
\]
Since $r_k=k+\alpha-1$ and $\sum_{k=2}^{+\infty}k\varepsilon_k<+\infty$, we have
$\sum_{k=1}^{+\infty}r_k\varepsilon_{k+1}<+\infty$. Lemma~\ref{le_disc_per} therefore gives $\sup_{k\geq1} \norm{
    U_{k+1}+\delta(x_{k+1}-x_k)
  }<+\infty$. Using Lemma~\ref{lem:dis_est} once more,
\[
  \sup_{k\geq1} {\mathcal E_{k+1,z^\star}}
  \leq
  {\mathcal E_{1,z^\star}}
  +
  \tau \sup_{k\geq1} \norm{
    U_{k+1}+\delta(x_{k+1}-x_k)
  }
  \sum_{k=1}^{+\infty}
  r_k\varepsilon_{k+1}<+\infty.
\]
Together with \eqref{eq:dis_energy}, this gives
$\sup_{k\geq1}\|U_{k+1}\|<+\infty$ and
$\sup_{k\geq1}\|V_{k+1}\|<+\infty$. Applying
Lemma~\ref{lem:bounded_dis} to the primal and dual sequences yields the
boundedness of $\{(x_k,\lambda_k)\}$ and
\[\|x_k-x_{k-1}\|+\|\lambda_k-\lambda_{k-1}\|=\mathcal O(k^{-1}).\]
This proves~(i).

(ii) We first establish the feasibility rate. Define, using \eqref{eq:dis}\textnormal{(e)},
\begin{equation}\label{eq_etadef}
  \eta_{k+1}
  :=
  \frac1\beta
  \bigl(
    p_{k+1}
    -\widehat\lambda_{k+1}
    -\beta\widehat g_{k+1}
  \bigr) =\frac1\beta
  \pos{
    -\widehat\lambda_{k+1}
    -\beta\widehat g_{k+1}
  }
  \in\R_+^m.
\end{equation}
Equation~\eqref{eq:dis}\textnormal{(b)} becomes
\begin{equation}\label{eq:feas_rel}
  \lambda_{k+1}-\bar\lambda_k
  =
  \tau\sigma
  \bigl(
    \widehat g_{k+1}+\eta_{k+1}
  \bigr).
\end{equation}

Define
\begin{equation}\label{eq_defG}
    G_k:=r_{k-1}^2g(x_k),
  \qquad
  \Lambda_k:=
  (\alpha-1)\lambda_k
  +(k-1)(\lambda_k-\lambda_{k-1}).
\end{equation}
By \eqref{eq:iner}, \eqref{eq:dis}, and \eqref{eq_rk}, we have
\[
  \Lambda_{k+1}-\Lambda_k=
  r_k(\lambda_{k+1}-\bar\lambda_k)\] and 
  \[
  \gamma\widehat g_{k+1}
  =
  r_k g(x_{k+1})
  -(r_k-\gamma)g(x_k).
\]
Combining these identities with \eqref{eq:feas_rel} gives
\[
  r_k^2g(x_{k+1})
  -r_k(r_k-\gamma)g(x_k)=\gamma r_k \widehat g_{k+1}
 =
  \frac{\gamma}{\tau\sigma}
  (\Lambda_{k+1}-\Lambda_k)
  -\gamma r_k\eta_{k+1}.
\]
Multiplying both sides by $
  \frac{r_{k-1}^2}{r_k(r_k-\gamma)}
$ and using \eqref{eq_rk_eqa} and \eqref{eq_defG}, we obtain
\begin{equation}\label{eq_diffG}
  G_{k+1}-G_k
  +
  \frac{(\gamma-2)r_k+1}
     {r_k(r_k-\gamma)}
  G_{k+1}
 =
  \frac{\gamma r_{k-1}^2}
     {\tau\sigma r_k(r_k-\gamma)}
  (\Lambda_{k+1}-\Lambda_k)
  -
  \frac{\gamma r_{k-1}^2}
     {r_k-\gamma}
  \eta_{k+1}.
\end{equation}
Define  $a_1:=0$ and 
\begin{equation}\label{eq_defab}
   \ 
  a_{k+1}:=
  \frac{(\gamma-2)r_k+1}{r_k(r_k-\gamma)},
  \qquad
  b_k:=
  \frac{\gamma r_{k-1}^2}{\tau\sigma r_k(r_k-\gamma)},\
  \forall k\geq1.
\end{equation}
By \eqref{eq_rk} and the parameter conditions, $a_k\geq0$ for all
$k\geq1$. Summing \eqref{eq_diffG}  from $k=1$ to $N-1$  gives, for every $N\geq1$,
\begin{align}\label{eq:dis1}
  G_N
  +\sum_{k=2}^{N}a_{k}G_{k}
  =
  G_1
  +\sum_{k=1}^{N-1}b_k(\Lambda_{k+1}-\Lambda_k)
-
  \sum_{k=1}^{N-1}
  \frac{\gamma r_{k-1}^2}{r_k-\gamma}\eta_{k+1}.
\end{align}
Moreover, $\{b_k\}$ is nonnegative and bounded, and
\begin{align}\label{eq_diffb}
  b_{k+1}-b_k
  &=
  \frac{\gamma}{\tau\sigma}
  \left[
    \frac{(k+\alpha-1)^2}
    {(k+\alpha)(k+\alpha-\gamma)}
    -
    \frac{(k+\alpha-2)^2}
    {(k+\alpha-1)(k+\alpha-1-\gamma)}
  \right]\nonumber \\
  &=-\frac{\gamma}{\tau\sigma}
  \frac{
    (\gamma-2)(k+\alpha-1)^2
    +\gamma(k+\alpha-1)-\gamma+1
  }{
    (k+\alpha-1)(k+\alpha-1-\gamma)
    (k+\alpha)(k+\alpha-\gamma)
  }.
\end{align}
Since $2\leq\gamma\leq\alpha-1$, there exists $C>0$ such that
$
  |b_{k+1}-b_k|\leq C/k^2
$
for every $k\geq1$. By part~(i), $\{\lambda_k\}$ is bounded and
$\|\lambda_k-\lambda_{k-1}\|=\mathcal O(k^{-1})$; hence
\eqref{eq_defG} implies that $\{\Lambda_k\}$ is bounded. For $N\geq2$,
summation by parts gives
\begin{align}\label{eq_bL1}
  \sum_{k=1}^{N-1}b_k(\Lambda_{k+1}-\Lambda_k)
  ={}&
  b_{N-1}\Lambda_N-b_1\Lambda_1
  +\sum_{k=2}^{N-1}(b_{k-1}-b_k)\Lambda_k.
\end{align}
The right-hand side is uniformly bounded because $\{b_k\}$ and
$\{\Lambda_k\}$ are bounded and
$\sum_{k=1}^{+\infty}|b_{k+1}-b_k|\le \sum_{k=1}^{+\infty} C/k^2 <+\infty$. Therefore, the sequence
\begin{align}\label{eq_bL2}
  B_N:=G_1+
  \sum_{k=1}^{N-1}b_k(\Lambda_{k+1}-\Lambda_k),
  \qquad N\geq1,
\end{align}
is bounded. On the other hand,
\[
  M_N:=
  \sum_{k=1}^{N-1}
  \frac{\gamma r_{k-1}^2}{r_k-\gamma}\eta_{k+1},
\]
satisfies $M_1=0$ and is componentwise nondecreasing, since
$r_k-\gamma>0$ and $\eta_{k+1}\in\R_+^m$. Then, from \eqref{eq:dis1} we have
\begin{align}\label{eq:disnew}
  G_N
  +\sum_{k=2}^{N}a_{k}G_{k}
  = B_N- M_N.
\end{align}
 Applying
Lemma~\ref{lem:dis_vol}(i) to \eqref{eq:disnew} therefore yields
$
  \sup_{k\geq1}\|\pos{G_k}\|<+\infty.
$
Consequently, by \eqref{eq_defG} and \eqref{eq_rk},
\begin{equation}\label{eq:dis_es4}
  \|\pos{g(x_k)}\|
  =\frac{\|\pos{G_k}\|}{r_{k-1}^2}
  =\mathcal O(k^{-2}).
\end{equation}

 We next estimate the objective residual. Keeping the fixed primal solution $x^\star$, we evaluate the energy expression \eqref{eq:dis_energy} at the reference pair $(x^\star,0)$ and define
\begin{equation}\label{eq_disEK0}
  \mathcal E_{k,0}
  :=
  \mathcal E_{k,z^\star}\big|_{\lambda^\star=0}
  =
  \tau r_{k-1}^2
  \bigl(f(x_k)-f^\star\bigr)
  +
  \mathcal Q_k,
\end{equation}
where
\[
\mathcal Q_k
=
\frac12\norm{U_k}^2
+\frac{\gamma\delta}{2}\norm{x_k-x^\star}^2
+\frac1{2\sigma}
\norm{
  \gamma\lambda_k
  +(k-1)(\lambda_k-\lambda_{k-1})
}^2
+\frac{\gamma\delta}{2\sigma}\norm{\lambda_k}^2.
\]
This is the energy \eqref{eq:dis_energy} with the 
multiplier $\lambda^\star$ replaced by $0$. Repeating the calculation in the proof of
Lemma~\ref{lem:dis_est} with the reference multiplier set to
$0\in\R_+^m$ gives the estimate below. This calculation uses feasibility
of $x^\star$ and the projection inequality with reference multiplier $0$;
it does not require $(x^\star,0)$ to be a KKT pair or
$f(x_k)-f^\star$ to be nonnegative. Thus
\begin{align}\label{eq:dis_obj}
    {\mathcal E_{k+1,0}-\mathcal E_{k,0}}
    \notag&\leq{}
    -\tau\bigl((\gamma-2)r_k+1\bigr)
    \bigl(f(x_k)-f^\star\bigr)
    \notag\\
    &\quad+\tau r_k
    \ip{\epsilon_{k+1}}
       {U_{k+1}+\delta(x_{k+1}-x_k)}.
\end{align}
By part~(i),
$
  0\leq\mathcal Q_k\leq C_0
$
for some $C_0>0$. Set
$
  \theta_k=
  \frac{(\gamma-2)r_k+1}{r_{k-1}^2}.
$
By \eqref{eq_rk}, \eqref{eq_rk_eqa}, and $\gamma\geq2$, we have
$0<\theta_k<1$ for every $k\geq1$. By \eqref{eq_disEK0},
\[
  f(x_k)-f^\star
  =
  \frac{{\mathcal E_{k,0}}-\mathcal Q_k}
     {\tau r_{k-1}^2}.
\]
Substituting this identity into
\eqref{eq:dis_obj} gives
\begin{align}\label{eq:E0_rec}
  {\mathcal E_{k+1,0}}
  &\leq
  (1-\theta_k){\mathcal E_{k,0}}
  +\theta_k\mathcal Q_k
  +
  \tau r_k
  \ip{\epsilon_{k+1}}
    {U_{k+1}+\delta(x_{k+1}-x_k)}
  \notag\\
  &\leq
   \max\{{\mathcal E_{k,0}},C_0\}
  + \tau\sup\nolimits_{j\geq1} \norm{
    U_{j+1}+\delta(x_{j+1}-x_j)
  }\cdot r_k\varepsilon_{k+1}\\
  &\leq
  \max\{{\mathcal E_{1,0}},C_0\}
  +
  \tau\sup\nolimits_{j\geq1} \norm{
    U_{j+1}+\delta(x_{j+1}-x_j)
  }\cdot\sum_{j=1}^{k}r_j\varepsilon_{j+1}.
  \nonumber
\end{align}
By (i) and the assumption,
\[
    \sup_{j\geq1}
    \|U_{j+1}+\delta(x_{j+1}-x_j)\|<+\infty,
\]
and $\sum_{k=1}^{+\infty}r_k\varepsilon_{k+1}<+\infty$.
Equation~\eqref{eq:E0_rec} therefore implies that
$
  {\mathcal E_{k,0}}
  \leq C_1
$
for some constant $C_1>0$. Hence, since $\mathcal Q_k\geq0$, \eqref{eq_disEK0} gives
\begin{equation}\label{eq:upp}
  f(x_k)-f^\star
  =
  \frac{{\mathcal E_{k,0}}-\mathcal Q_k}
     {\tau r_{k-1}^2}
  \leq
  \frac{C_1}{\tau r_{k-1}^2}.
\end{equation}

On the other hand, the saddle-point inequality \eqref{eq:saddle} gives
\[
  f(x_k)-f^\star
  +\ip{\lambda^\star}{g(x_k)}
  \geq0.
\]
Since $\lambda^\star\in\R_+^m$,
\[
\begin{aligned}
  f(x_k)-f^\star
  &\geq
  -\ip{\lambda^\star}{g(x_k)}
  \geq
  -\ip{\lambda^\star}{\pos{g(x_k)}}\geq
  -\norm{\lambda^\star}
   \norm{\pos{g(x_k)}}.
\end{aligned}
\]
Combining this estimate with
\eqref{eq:dis_es4} and \eqref{eq:upp} yields
$
  |f(x_k)-f^\star|
  =
  \mathcal O(k^{-2}).
$
This proves~(ii).

\end{proof}

\begin{remark}
Theorem~\ref{thm:dis_rates} gives simultaneous $\mathcal O(k^{-2})$ absolute objective and feasibility estimates at the generated primal iterates, including $\alpha=3$. It requires neither the objective strong convexity used for acceleration in \cite{HamedaniAybat2021,ConEx2023,ZhuSiopt} nor the constraint strong convexity used in APDPro \cite{LinNeurips}. In the general convex regime,  ACGD \cite{ZhangLanACGD} gives weighted-average estimates, and smooth APMM \cite{DengInformsJoc} gives a merit bound at its prescribed primal output. Here both rates hold at the generated iterates, allow a nonsmooth objective component, and accommodate weighted-summable primal stationarity errors.
\end{remark}

\subsection{Convergence of the iterates and improved rates}
 
 Under strict parameter inequalities, the energy dissipation yields convergence of the  primal-dual sequence and little-$o$ residual rates.

Throughout this part, in addition to the assumptions of
Theorem~\ref{thm:dis_rates}, suppose that
\begin{equation}\label{eq:strict_para}
  \alpha>3,\qquad 2<\gamma<\alpha-1.
\end{equation}
Thus $\delta=\alpha-\gamma-1>0$. We use the notation
\[
  z_k=(x_k,\lambda_k),\qquad
  \|z\|_\sigma^2=\|x\|^2+\frac1\sigma\|\lambda\|^2
  \quad\text{for }z=(x,\lambda),
\]
and \[d_{k+1}:=p_{k+1}-\widehat\lambda_{k+1}.\]
For $z^\star=(x^\star,\lambda^\star)\in\Omega$, we continue to use
the notation $\mathcal E_{k,z^\star}$ from \eqref{eq:dis_energy}, with
$\ell_k,U_k,V_k$ in \eqref{eq:aux} computed relative to this pair.

\begin{lemma}\label{lem:dis_cluster}
Suppose that the assumptions of Theorem~\ref{thm:dis_rates} and
\eqref{eq:strict_para} hold. Then every cluster point of
$\{(x_k,\lambda_k)\}$ generated by Algorithm~\ref{alg:inexact}
belongs to $\Omega$.
\end{lemma}

\begin{proof}
Fix $z^\star=(x^\star,\lambda^\star)\in\Omega$.
By the proof of Theorem~\ref{thm:dis_rates}(i), there exists $C>0$ such that
$\sup_{k\ge 1}\|U_{k+1}+\delta(x_{k+1}-x_k)\|\le C$. Hence the
inexactness assumption gives
\[  \sum_{k=1}^{+\infty}
  \tau r_k
  \left|
    \ip{\epsilon_{k+1}}
      {U_{k+1}+\delta(x_{k+1}-x_k)}
  \right|
  \leq
  \tau C\sum_{k=1}^{+\infty}r_k\varepsilon_{k+1}
  <+\infty.
\]
Since
$\gamma>2$, $\delta>0$, and $r_k\geq k$, summing \eqref{eq:dis_full_est} over $k\geq1$ and using ${\mathcal E_{k,z^\star}}\geq0$, we obtain
\begin{equation}\label{eq_sumres}
  \sum_{k=1}^{+\infty}
  k\left(
    \ell_k+\|z_{k+1}-z_k\|_\sigma^2
    +\|d_{k+1}\|^2
  \right)<+\infty.
\end{equation}

Let $\bar z=(\bar x,\bar\lambda)$ be an arbitrary cluster
point of $\{z_k\}$, and choose $n_j\to+\infty$ such that
$z_{n_j}\to\bar z$. Fix $\rho>1$. By the preceding
summability and the equivalence of norms in finite dimensions,
Lemma~\ref{lem:est_diff} gives
\[  \sup_{n_j\leq q\leq\lfloor\rho n_j\rfloor}
  \|z_q-\bar z\|\to0
\]
and there exist integers $q_j\in[n_j,\lfloor\rho n_j\rfloor]$ such that
\begin{equation}\label{eq:dis_ind}
  q_j\left(
    \|z_{q_j}-z_{q_j-1}\|
    +\|z_{q_j+1}-z_{q_j}\|
    +\|d_{q_j+1}\|
  \right)\to0.
\end{equation}
This, together with $z_{n_j}\to\bar z$, implies
\[
  z_{q_j}\to\bar z,\qquad
  z_{q_j+1}
  =z_{q_j}+(z_{q_j+1}-z_{q_j})\to\bar z.
\]
Since $\frac{r_{q_j}}{q_j}=\frac{q_j+\alpha-1}{q_j}\to1$ as $j\to+\infty$, \eqref{eq:dis}\textnormal{(c)}
and \eqref{eq:dis_ind} imply
\begin{equation}\label{eq_hatLto}
	  \widehat\lambda_{q_j+1}
  =
  \lambda_{q_j}
  +\frac{r_{q_j}}{\gamma}
   (\lambda_{q_j+1}-\lambda_{q_j})
  \to\bar\lambda,
\end{equation}
as $j\to+\infty$. By Theorem~\ref{thm:dis_rates}(i), $\{x_k\}$ is bounded.
Since $g$ is continuously differentiable, it is Lipschitz
continuous on a closed ball containing all the iterates.
Thus, by \eqref{eq:dis}\textnormal{(d)}, there exists
$C>0$, independent of $j$, such that
\[
  \|\widehat g_{q_j+1}-g(\bar x)\|
  \leq
  \|g(x_{q_j})-g(\bar x)\|
  +C r_{q_j}\|x_{q_j+1}-x_{q_j}\|.
\]
Together with
$x_{q_j}\to\bar x$, $r_{q_j}/q_j\to1$, and
\eqref{eq:dis_ind}, this estimate gives
\[\widehat g_{q_j+1}\to g(\bar x)\qquad\text{as }j\to+\infty.\]

Moreover, \eqref{eq:dis_ind} gives
$d_{q_j+1}\to0$, and hence, by \eqref{eq_hatLto},
\begin{equation}\label{eq_pklim}
	  p_{q_j+1}
  =\widehat\lambda_{q_j+1}+d_{q_j+1}
  \to\bar\lambda.
\end{equation}
Passing to the limit in \eqref{eq:dis}\textnormal{(e)}
along the subsequence $\{q_j\}$ as $j\to+\infty$ therefore yields
\begin{equation}\label{eq:clusterP}
  \bar\lambda=\pos{\bar\lambda+\beta g(\bar x)}.
\end{equation}
This is equivalent to
\[
  g(\bar x)\leq0,\qquad
  \bar\lambda\geq0,\qquad
  \ip{\bar\lambda}{g(\bar x)}=0.
\]

It remains to verify stationarity.
By \eqref{eq:iner} and \eqref{eq:dis_ind},
\[
  x_{q_j+1}-\bar x_{q_j}
  =
  x_{q_j+1}-x_{q_j}
  -\frac{q_j-1}{r_{q_j}}(x_{q_j}-x_{q_j-1})
  \to0.
\]
Together with $x_{q_j+1}\to\bar x$, this gives
$\bar x_{q_j}\to\bar x$. Choose
$\xi_{q_j+1}\in\partial h(x_{q_j+1})$ as in
\eqref{eq:res}. Since
$\|\epsilon_{q_j+1}\|\leq\varepsilon_{q_j+1}\to0$,
the continuity of $\nabla\phi$ and $J_g$ together with \eqref{eq_pklim} gives
\begin{align*}
  \xi_{q_j+1}
  &=
  \epsilon_{q_j+1}
  -\frac{x_{q_j+1}-\bar x_{q_j}}{\tau}
  -\nabla\phi(\bar x_{q_j})
  -J_g(x_{q_j+1})^\top p_{q_j+1}\\
  &\to
  -\nabla\phi(\bar x)-J_g(\bar x)^\top\bar\lambda.
\end{align*}
as $j\to+\infty$. The subdifferential of a proper closed convex function has a closed
graph in finite dimensions \cite{Bauschke16}.
Since
$x_{q_j+1}\to\bar x$ as $j\to+\infty$, it follows that
\[
  0\in
  \nabla\phi(\bar x)+\partial h(\bar x)
  +J_g(\bar x)^\top\bar\lambda.
\]
Together with \eqref{eq:clusterP} and \eqref{eq:KKT}, this proves
$\bar z\in\Omega$.
\end{proof}

\begin{theorem}\label{thm:dis_littleo}
Let $\{(x_k,\lambda_k)\}$ be generated by Algorithm~\ref{alg:inexact}.
Suppose that the assumptions of Theorem~\ref{thm:dis_rates} and
\eqref{eq:strict_para} hold. Then the following statements hold.
\begin{itemize}
  \item[(i)] There exists
  $(x^\infty,\lambda^\infty)\in\Omega$ such that
  \[
    (x_k,\lambda_k)\to(x^\infty,\lambda^\infty)
    \qquad\text{as }k\to+\infty.
  \]
  \item[(ii)] The iterates satisfy
  \[
    \|x_k-x_{k-1}\|+\|\lambda_k-\lambda_{k-1}\|
    =o\left(\frac{1}{k}\right).
  \]
  \item[(iii)] The feasibility and objective residuals satisfy
\[ \|\pos{g(x_k)}\|=o\left(\frac{1}{k^{2}}\right),\qquad |f(x_k)-f^\star|=o\left(\frac{1}{k^{2}}\right). \]
\end{itemize}\end{theorem}

\begin{proof}
(i) For any $z^\star=(x^\star,\lambda^\star)\in\Omega$, set
\[
  h_k^{z^\star}
  :=
  \frac12\|z_k-z^\star\|_\sigma^2
\]
and
\begin{equation}\label{eq_defR}
	  R_k^{z^\star}
  :=
  \tau r_{k-1}^2\ell_k
  +
  \frac{(k-1)(k+\gamma-1)}2
  \|z_k-z_{k-1}\|_\sigma^2,
\end{equation}
where $\ell_k$ is computed with the same KKT pair. Expanding the two
squares in \eqref{eq:dis_energy}, and using
$\gamma+\delta=\alpha-1$, gives the exact identity
\begin{equation}\label{eq:energy_dis}
  {\mathcal E_{k,z^\star}}
  =
  R_k^{z^\star}
  +\gamma(k-1)
   \bigl(h_k^{z^\star}-h_{k-1}^{z^\star}\bigr)
  +\gamma(\alpha-1)h_k^{z^\star}.
\end{equation}
The first two summability estimates in \eqref{eq_sumres}
imply
\begin{equation}\label{eq:dis_R_sum}
  \sum_{k=2}^{+\infty}\frac{R_k^{z^\star}}{k}<+\infty.
\end{equation}
Rearranging \eqref{eq:energy_dis} gives
\begin{equation}\label{eq:diffhk}
  (k-1)
  \bigl(h_k^{z^\star}-h_{k-1}^{z^\star}\bigr)
  +(\alpha-1)h_k^{z^\star}
  =
  \frac1\gamma{\mathcal E_{k,z^\star}}
  -\frac1\gamma R_k^{z^\star}.
\end{equation}

By Theorem~\ref{thm:dis_rates}(i), there exists $C>0$ such that
\[
  \|U_{k+1}+\delta(x_{k+1}-x_k)\|\leq C
  \qquad\text{for all }k\geq1.
\]
Since all dissipative terms in \eqref{eq:dis_full_est}
are nonpositive, we obtain
\[
  {\mathcal E_{k+1,z^\star}-\mathcal E_{k,z^\star}}
  \leq
  \tau r_k
  \ip{\epsilon_{k+1}}
     {U_{k+1}+\delta(x_{k+1}-x_k)}
  \leq C\tau r_k\varepsilon_{k+1}.
\]
The inexactness assumption gives $\sum_{k=1}^{+\infty}r_k\varepsilon_{k+1}<+\infty$. The preceding inequality shows that $\mathcal E_{k,z^\star}+C\tau\sum_{j=k}^{+\infty}r_j\varepsilon_{j+1}$ is nonnegative and nonincreasing, so it has a finite limit. Since the summable tail tends to zero, $\mathcal E_{k,z^\star}$ also has a finite limit. Together with
\eqref{eq:dis_R_sum}, Lemma~\ref{lem:w_limit}, applied to \eqref{eq:diffhk}
with $a=\alpha-1$, $E_k={\mathcal E_{k,z^\star}}/\gamma$, and
$R_k=R_k^{z^\star}/\gamma$, shows that, for every $z^\star\in\Omega$,
\begin{equation}\label{eq:dis_limit}
  \lim_{k\to+\infty}h_k^{z^\star}
  =
  \frac{
    \lim_{k\to+\infty}{\mathcal E_{k,z^\star}}
  }{\gamma(\alpha-1)}\quad \text{exists}.
\end{equation}
In particular, $\|z_k-z^\star\|_\sigma$ has a limit.

By Lemma~\ref{lem:dis_cluster}, every cluster point belongs to $\Omega$.
The discrete version of Opial's lemma (Lemma~\ref{lem:opial_dis}), applied
with the fixed norm $\|\cdot\|_\sigma$, therefore gives a pair
$z^\infty=(x^\infty,\lambda^\infty)\in\Omega$ such that
  \[
      z_k\to z^\infty
    \qquad\text{as }k\to+\infty.
  \]
This proves~(i).
   
(ii)  Take $z^\star=z^\infty$. Then $h_k^{z^\infty}\to0$, and 
equation~\eqref{eq:dis_limit} with
$z^\star=z^\infty$ also gives
\begin{equation}\label{eq:lim_energy}
  \lim_{k\to+\infty}{\mathcal E_{k,z^\infty}}=0.
\end{equation}
Directly from the definition of
$h_k^{z^\infty}$,
\begin{align*}
 (k-1)
 \bigl(h_k^{z^\infty}-h_{k-1}^{z^\infty}\bigr) &= (k-1)\ip{x_k-x^\infty}{x_k-x_{k-1}}
 +\frac{k-1}{\sigma}
  \ip{\lambda_k-\lambda^\infty}{\lambda_k-\lambda_{k-1}}\\
 &\quad
 -\frac{k-1}{2}
 \left(
   \|x_k-x_{k-1}\|^2
   +\frac1\sigma\|\lambda_k-\lambda_{k-1}\|^2
 \right).
\end{align*}
Together with Theorem~\ref{thm:dis_rates}(i) and $z_k\to z^\infty$, this proves
\[
  (k-1)
  \bigl(h_k^{z^\infty}-h_{k-1}^{z^\infty}\bigr)
  \to0.
\]
Using this limit, $h_k^{z^\infty}\to0$, and
\eqref{eq:lim_energy} and \eqref{eq:energy_dis}, we obtain
\[
 \lim_{k\to+\infty} R_k^{z^\infty}=0.
\]
It follows from \eqref{eq_defR} and the strict parameter inequalities that
\begin{equation}\label{eq:dis_gapo}
  \mathcal L(x_k,\lambda^\infty)-f^\star=o(k^{-2}),
  \qquad
  \|x_k-x_{k-1}\|+\|\lambda_k-\lambda_{k-1}\|
  =o(k^{-1}).
\end{equation}
This proves (ii).

(iii) Use the quantities $G_k,\Lambda_k,a_k,b_k,B_k,M_k$, and $\eta_k$ introduced
in the proof of Theorem~\ref{thm:dis_rates}(ii). It follows from \eqref{eq_defG}, (i) and (ii) that
\begin{equation}\label{eq_Llim}
  \Lambda_k
  =
  (\alpha-1)\lambda_k
  +(k-1)(\lambda_k-\lambda_{k-1})
  \to(\alpha-1)\lambda^\infty.
\end{equation}
Moreover, \eqref{eq_defab} and \eqref{eq_diffb} give
\begin{equation}\label{eq_blim}
 \lim_{k\to+\infty} b_k = \frac{\gamma}{\tau\sigma},
  \qquad
  \sum_{k=1}^{+\infty}|b_{k+1}-b_k|<+\infty.
\end{equation}
By \eqref{eq_bL1} and \eqref{eq_bL2}, for $k\geq2$,
\[
  B_k
  =
  G_1+b_{k-1}\Lambda_k-b_1\Lambda_1
  +\sum_{j=2}^{k-1}(b_{j-1}-b_j)\Lambda_j.
\]
It follows from \eqref{eq_Llim} and \eqref{eq_blim} that
$\{\Lambda_k\}$ is bounded. Hence
\[
  \sum_{j=2}^{+\infty}
  \|(b_{j-1}-b_j)\Lambda_j\|
  \leq
  \left(\sup_{j\geq1}\|\Lambda_j\|\right)
  \sum_{j=2}^{+\infty}|b_{j-1}-b_j|
  <+\infty.
\]
 Moreover,
\[
  \lim_{k\to+\infty}b_{k-1}\Lambda_k
  =
  \frac{\gamma(\alpha-1)}{\tau\sigma}\lambda^\infty.
\]
Since $G_1$ and $b_1\Lambda_1$ are fixed, it follows
that $B_k$ has a finite limit. In addition,
\[
  a_{k+1}
  =
  \frac{(\gamma-2)r_k+1}{r_k(r_k-\gamma)}
  \sim\frac{\gamma-2}{k},
\]
so $\sum_{k=2}^{+\infty}a_k=+\infty$. Recall from the proof of Theorem~\ref{thm:dis_rates}(ii) that
\[
  G_k+\sum_{j=2}^k a_jG_j=B_k-M_k
\]
with $M_k
  :=
  \sum_{j=1}^{k-1}
  \frac{\gamma r_{j-1}^2}{r_j-\gamma}\eta_{j+1}.$
As verified in that proof, $a_k\geq0$, $M_1=0$, and
$\{M_k\}$ is componentwise nondecreasing.
Together with the finite limit of $\{B_k\}$ and
$\sum_{k=2}^{+\infty}a_k=+\infty$ established above,
these properties allow us to apply
Lemma~\ref{lem:dis_vol}(ii), which gives
\[
  \lim_{k\to+\infty}\|\pos{G_k}\|=0.
\]
Since $G_k=r_{k-1}^2g(x_k)$ and the positive-part mapping is positively
homogeneous,
\[  \|\pos{g(x_k)}\|=o\left(\frac{1}{k^{2}}\right).
\]

Finally, we prove the objective residual estimate.
Since $r_k/(k+1)\to1$, (ii) gives $r_k\|x_{k+1}-x_k\|\to0$ and $r_k\|\lambda_{k+1}-\lambda_k\|\to0$. It follows from (i) and \eqref{eq:dis}\textnormal{(c)} that
\[
  \widehat\lambda_{k+1}
  =\lambda_k+\frac{r_k}{\gamma}
  (\lambda_{k+1}-\lambda_k)
  \to\lambda^\infty.
\]
By (i), a closed ball $B$ contains all the primal iterates and $x^\infty$. Since $g$ is continuously differentiable, compactness of $B$ gives $L_B:=\max_{x\in B}\|J_g(x)\|<+\infty$. Integrating the Jacobian along the segment joining any $u,v\in B$ yields $\|g(u)-g(v)\|\leq L_B\|u-v\|$. Combining the bound with \eqref{eq:dis}\textnormal{(d)} and the continuity of $g$ yields
\[
  \|\widehat g_{k+1}-g(x^\infty)\|
  \leq
  \|g(x_k)-g(x^\infty)\|
  +\frac{L_Br_k}{\gamma}\|x_{k+1}-x_k\|
  \to0.
\]

If $\lambda_i^\infty>0$, complementarity gives
$g_i(x^\infty)=0$, and hence
\[
  \widehat\lambda_{k+1,i}
  +\beta\widehat g_{k+1,i}
  \to\lambda_i^\infty>0.
\]
Thus $\eta_{k+1,i}=0$ for all sufficiently large $k$, where $\eta_k$ is defined in \eqref{eq_etadef}. Consequently,
there exists an integer $K_i$ such that
$M_{k,i}=M_{K_i,i}$ for every $k\geq K_i$. The last conclusion of
Lemma~\ref{lem:dis_vol}(ii) gives
\[
  G_{k,i}\to0
  \qquad\text{whenever }\lambda_i^\infty>0.
\]
If $\lambda_i^\infty=0$, then $\lambda_i^\infty G_{k,i}=0$ for every $k$. Since $\lambda^\infty\geq0$, summing the finitely many remaining terms gives $\ip{\lambda^\infty}{G_k}=\sum_{i:\,\lambda_i^\infty>0}\lambda_i^\infty G_{k,i}\to0$.
Finally, by \eqref{eq:dis_gapo} and $G_k=r_{k-1}^2g(x_k)$,
\begin{align*}
  r_{k-1}^2\bigl(f(x_k)-f^\star\bigr)
  &=
  r_{k-1}^2
  \bigl(\mathcal L(x_k,\lambda^\infty)-f^\star\bigr)
  -\ip{\lambda^\infty}{G_k}\to0.
\end{align*}
This proves $|f(x_k)-f^\star|=o(k^{-2})$ and  completes the proof.
\end{proof}

\begin{remark}
 For $\alpha>3$ and $2<\gamma<\alpha-1$, Theorem~\ref{thm:dis_littleo} gives simultaneous last-iterate $o(k^{-2})$ objective and feasibility rates, convergence of the entire primal-dual sequence, and $o(k^{-1})$ successive differences. In the general convex regime, ACGD \cite{ZhangLanACGD} gives $\mathcal O(k^{-2})$ rates for weighted primal averages, and smooth APMM \cite{DengInformsJoc} gives an $\mathcal O(k^{-2})$ merit bound. The same weighted error condition suffices, with no strong convexity or uniqueness assumption. This extends little-$o$ acceleration and iterate convergence from unconstrained \cite{AttouchPtSiopt,ChambolleJota} and affine-constrained \cite{HeArxiv} problems to nonlinear convex inequalities. When $g\equiv0$ and primal solves are exact, the primal step reduces to \eqref{eq:fista} and recovers the corresponding unconstrained guarantees \cite{AttouchMathProg,AttouchPtSiopt,ChambolleJota}.
\end{remark}

\section{Numerical experiments}\label{sec:numerics}
 In this section, we conduct numerical experiments to evaluate the effectiveness of the proposed methods. All implementations use MATLAB 9.4.0.813654 (R2018a) on a PC with a 2.3 GHz quad-core Intel Core i5 processor and 8 GB of memory.

 \subsection{Weighted sparse recovery}
 We use a weighted variant of the basis-pursuit denoising model studied in \cite{Berg2011}:
\[ \min_{x\in\R^n}\ f(x):=\sum_{i=1}^n w_i|x_i|
 \quad\text{subject to}\quad \|Ax-b\|_2\leq\delta.
\]
Here $A\in\mathbb R^{m\times n}$. We take $n=3000$ and $m=1000$.
The entries of $A$ are independent standard Gaussian random variables before
rescaling to $\|A\|_2=1$, and the weights are independent and uniformly
distributed on $[0.7,1.3]$. The vector $x^{\rm true}$ has $0.1m$ randomly
located nonzero standard Gaussian entries. We set $b=Ax^{\rm true}+\xi$,
where $\xi_i\sim N(0,0.004^2)$, and $\delta=\|\xi\|$.
All methods use the same instance and start from zero primal and
dual variables. In the notation of the present paper, $\phi=0$ and $h=f$;
the primal domain is $\mathbb R^n$. We replace the constraint by the
equivalent smooth inequality
\[
 g(x):=\frac{\|Ax-b\|^2-\delta^2}{2s_b}\leq0,
 \qquad s_b=\max\{\|b\|,\delta\}.
\]
An independent primal-dual conic solve provides primal- and dual-feasible
reference points $x^{\rm ref}$ and $y^{\rm ref}$, from which $f_{\rm ref}$
is computed with absolute accuracy $10^{-12}$.

We compare AALM with the inexact augmented Lagrangian method \cite[Algo.~1]{XuMathProg} (iALM), the adaptively regularized inexact augmented Lagrangian method \cite[Algo.~2]{LuSiopt} (RI-AL), and the linearized augmented Lagrangian method \cite[Algo.~1]{XuLinear2021} (LALM). For AALM, we use
$
 (\alpha,\gamma,\beta,\sigma,\tau)=(m,m-2,0.01,10m,100).
$
Here $L_\phi=0$, so the choice $\tau=100$ is admissible. The weighted $\ell_1$ objective is not strongly convex; neither is the quadratic constraint, since $m<n$. AALM treats the nonlinear penalty implicitly, and its outer stepsize bound involves only $L_\phi$.
For iALM, we set $\beta_k=20000$ and $\rho_k=20000$.
For RI-AL, we set $\rho_k=100(1.1)^k$.
For LALM, we set $\beta=200$ and $\rho_z=200$.
The strongly convex subproblems of AALM and RI-AL are solved using a
strongly convex accelerated proximal gradient method, whereas the convex subproblems of iALM are solved by the standard fast iterative shrinkage-thresholding algorithm (FISTA) \cite{BeckSiis,NesterovBook}.
LALM's subproblem has a closed-form proximal solution and requires no
iterative inner solve. For every iterative subproblem solve, we use the
same fixed conditions
\[ 
\operatorname{dist}(0,\partial F_k(x_{k+1}))\leq e_k,
 \qquad\text{at most }10^6\text{ inner updates},
\]
where $F_k$ is the objective function of the corresponding algorithm's
subproblem at iteration $k$. LALM uses its closed-form proximal map and
has zero iterative inner updates. For each fixed inner tolerance $e_k \in\{10^{-3},10^{-5},10^{-7}\}$,
we report the first outer iterate satisfying
\[
 R_K:=|f(x_K)-f_{\rm ref}|+[g(x_K)]_+\leq\varepsilon_{\rm out},
 \qquad\varepsilon_{\rm out}\in\{10^{-4},10^{-6},10^{-8}\}.
\]
 
\begin{table*}[t]
\centering
\caption{Weighted sparse recovery with $e_k=10^{-3}$.}
\label{tab:weighted:inner1}
\resizebox{\textwidth}{!}{%
\begin{tabular}{clrrrrrr}
\toprule
$\varepsilon_{\rm out}$ & Method & $K$ & $|f(x_K)-f_{\rm ref}|$ & $[g(x_K)]_+$ & $\|x_K-x_{K-1}\|$ & $\|\lambda_K-\lambda_{K-1}\|$ & Time (s) \\
\midrule
$10^{-4}$ & AALM & \textbf{11} & \texttt{5.18e-5} & \texttt{3.84e-8} & \texttt{5.01e-6} & \texttt{3.87e-2} & \textbf{1.799} \\
 & iALM & 23 & \texttt{3.64e-5} & \texttt{2.70e-8} & \texttt{1.05e-5} & \texttt{7.52e-3} & 1.880 \\
 & RI-AL & 64 & \texttt{6.88e-5} & \texttt{5.11e-8} & \texttt{9.77e-6} & \texttt{2.89e-2} & 2.115 \\
 & LALM & 3395 & \texttt{9.99e-5} & \texttt{0.00e0} & \texttt{1.43e-6} & \texttt{2.06e-4} & 5.787 \\
\midrule
$10^{-6}$ & AALM & \textbf{22} & \texttt{8.90e-7} & \texttt{6.61e-10} & \texttt{1.51e-6} & \texttt{6.72e-4} & \textbf{1.847} \\
 & iALM & 38 & \texttt{3.12e-7} & \texttt{0.00e0} & \texttt{2.92e-7} & \texttt{6.46e-5} & 1.936 \\
 & RI-AL & 77 & \texttt{4.17e-7} & \texttt{3.10e-10} & \texttt{3.37e-7} & \texttt{6.04e-4} & 2.170 \\
 & LALM & 4264 & \texttt{9.99e-7} & \texttt{0.00e0} & \texttt{3.02e-8} & \texttt{2.07e-6} & 7.215 \\
\midrule
$10^{-8}$ & AALM & \textbf{39} & \texttt{3.67e-9} & \texttt{0.00e0} & \texttt{1.61e-7} & \texttt{3.45e-6} & \textbf{1.921} \\
 & iALM & 50 & \texttt{5.82e-9} & \texttt{4.32e-12} & \texttt{3.23e-8} & \texttt{1.20e-6} & 1.981 \\
 & RI-AL & 110 & \texttt{8.66e-9} & \texttt{6.49e-12} & \texttt{8.69e-9} & \texttt{2.94e-4} & 2.306 \\
 & LALM & 5449 & \texttt{9.96e-9} & \texttt{0.00e0} & \texttt{1.74e-10} & \texttt{2.06e-8} & 9.195 \\
\bottomrule
\end{tabular}}
\end{table*}

\begin{table*}[t]
\centering
\caption{Weighted sparse recovery with $e_k=10^{-5}$.}
\label{tab:weighted:inner2}
\resizebox{\textwidth}{!}{%
\begin{tabular}{clrrrrrr}
\toprule
$\varepsilon_{\rm out}$ & Method & $K$ & $|f(x_K)-f_{\rm ref}|$ & $[g(x_K)]_+$ & $\|x_K-x_{K-1}\|$ & $\|\lambda_K-\lambda_{K-1}\|$ & Time (s) \\
\midrule
$10^{-4}$ & AALM & \textbf{11} & \texttt{4.04e-5} & \texttt{3.00e-8} & \texttt{5.89e-6} & \texttt{3.02e-2} & \textbf{2.384} \\
 & iALM & 24 & \texttt{6.63e-5} & \texttt{4.92e-8} & \texttt{2.60e-6} & \texttt{1.37e-2} & 2.944 \\
 & RI-AL & 64 & \texttt{9.09e-5} & \texttt{6.75e-8} & \texttt{7.27e-6} & \texttt{3.81e-2} & 3.939 \\
 & LALM & 3395 & \texttt{9.99e-5} & \texttt{0.00e0} & \texttt{1.43e-6} & \texttt{2.06e-4} & 5.787 \\
\midrule
$10^{-6}$ & AALM & \textbf{15} & \texttt{7.73e-7} & \texttt{5.74e-10} & \texttt{6.46e-8} & \texttt{5.81e-4} & \textbf{2.422} \\
 & iALM & 34 & \texttt{2.37e-7} & \texttt{1.76e-10} & \texttt{1.18e-7} & \texttt{4.90e-5} & 3.022 \\
 & RI-AL & 69 & \texttt{9.17e-7} & \texttt{6.81e-10} & \texttt{1.60e-7} & \texttt{6.19e-4} & 4.003 \\
 & LALM & 4264 & \texttt{9.99e-7} & \texttt{0.00e0} & \texttt{3.02e-8} & \texttt{2.07e-6} & 7.215 \\
\midrule
$10^{-8}$ & AALM & \textbf{25} & \texttt{6.19e-9} & \texttt{0.00e0} & \texttt{1.42e-8} & \texttt{4.50e-6} & \textbf{2.467} \\
 & iALM & 41 & \texttt{6.31e-10} & \texttt{0.00e0} & \texttt{2.45e-8} & \texttt{1.30e-7} & 3.049 \\
 & RI-AL & 77 & \texttt{1.07e-9} & \texttt{0.00e0} & \texttt{1.34e-8} & \texttt{1.55e-6} & 4.036 \\
 & LALM & 5449 & \texttt{9.96e-9} & \texttt{0.00e0} & \texttt{1.74e-10} & \texttt{2.06e-8} & 9.195 \\
\bottomrule
\end{tabular}}
\end{table*}

\begin{table*}[t]
\centering
\caption{Weighted sparse recovery with $e_k=10^{-7}$.}
\label{tab:weighted:inner3}
\resizebox{\textwidth}{!}{%
\begin{tabular}{clrrrrrr}
\toprule
$\varepsilon_{\rm out}$ & Method & $K$ & $|f(x_K)-f_{\rm ref}|$ & $[g(x_K)]_+$ & $\|x_K-x_{K-1}\|$ & $\|\lambda_K-\lambda_{K-1}\|$ & Time (s) \\
\midrule
$10^{-4}$ & AALM & \textbf{11} & \texttt{4.07e-5} & \texttt{3.02e-8} & \texttt{5.69e-6} & \texttt{3.04e-2} & \textbf{2.811} \\
 & iALM & 24 & \texttt{6.58e-5} & \texttt{4.88e-8} & \texttt{2.60e-6} & \texttt{1.36e-2} & 4.093 \\
 & RI-AL & 64 & \texttt{9.06e-5} & \texttt{6.72e-8} & \texttt{7.25e-6} & \texttt{3.80e-2} & 5.807 \\
 & LALM & 3395 & \texttt{9.99e-5} & \texttt{0.00e0} & \texttt{1.43e-6} & \texttt{2.06e-4} & 5.787 \\
\midrule
$10^{-6}$ & AALM & \textbf{15} & \texttt{5.84e-7} & \texttt{4.33e-10} & \texttt{8.12e-8} & \texttt{4.38e-4} & \textbf{2.983} \\
 & iALM & 34 & \texttt{9.05e-7} & \texttt{6.72e-10} & \texttt{3.54e-8} & \texttt{1.87e-4} & 4.513 \\
 & RI-AL & 70 & \texttt{3.50e-7} & \texttt{2.60e-10} & \texttt{5.05e-8} & \texttt{2.60e-4} & 6.070 \\
 & LALM & 4264 & \texttt{9.99e-7} & \texttt{0.00e0} & \texttt{3.02e-8} & \texttt{2.07e-6} & 7.215 \\
\midrule
$10^{-8}$ & AALM & \textbf{19} & \texttt{9.08e-9} & \texttt{6.74e-12} & \texttt{1.52e-9} & \texttt{6.83e-6} & \textbf{3.025} \\
 & iALM & 44 & \texttt{2.35e-9} & \texttt{1.75e-12} & \texttt{1.20e-9} & \texttt{4.87e-7} & 4.597 \\
 & RI-AL & 73 & \texttt{7.03e-9} & \texttt{5.22e-12} & \texttt{2.30e-9} & \texttt{6.95e-6} & 6.107 \\
 & LALM & 5449 & \texttt{9.96e-9} & \texttt{0.00e0} & \texttt{1.74e-10} & \texttt{2.06e-8} & 9.195 \\
\bottomrule
\end{tabular}}
\end{table*}

Tables~\ref{tab:weighted:inner1}--\ref{tab:weighted:inner3} show that AALM uses the fewest outer iterations and the shortest reported runtime at every prescribed accuracy. Tighter inner solves reduce the outer count at high accuracy, but their added cost increases the total runtime in this experiment. AALM already reaches the strictest outer tolerance with the loosest tested inner tolerance. LALM's explicit proximal update is independent of the inner tolerance, so the same run is reported in all three tables; its lower per-iteration cost does not offset the larger outer count on this instance.

\begin{table*}[t]\centering
\caption{Composite regression with $\mu=0.05$.}
\label{tab:comp1}
\resizebox{\textwidth}{!}{%
\begin{tabular}{clcrrrrr}\toprule
$\varepsilon_{\rm out}$ & Method & $e_k$ & $K$ & $|f(x_K)-f_{\rm ref}|$ & $\|[Cx_K-d]_+\|$ & $\|\Delta x_K\|+\|\Delta\lambda_K\|$ & Time (s) \\
\midrule
$10^{-4}$ & AALM & $10^{-3}$ & 33 & \texttt{6.51e-5} & \texttt{1.41e-5} & \texttt{4.71e-3} & \textbf{0.184} \\
 &  & $10^{-5}$ & \textbf{32} & \texttt{6.30e-5} & \texttt{1.40e-5} & \texttt{5.75e-3} & 0.285 \\
 &  & $10^{-7}$ & \textbf{32} & \texttt{6.27e-5} & \texttt{1.40e-5} & \texttt{5.75e-3} & 0.445 \\
\cmidrule(lr){2-8}
 & ConEx & -- & 722 & \texttt{9.98e-5} & \texttt{0.00e0} & \texttt{4.12e-5} & 2.321 \\
\cmidrule(lr){2-8}
 & AGRPDA-L & $10^{-3}$ & 143 & \texttt{9.70e-5} & \texttt{1.95e-6} & \texttt{1.36e-3} & 3.835 \\
 &  & $10^{-5}$ & 141 & \texttt{9.34e-5} & \texttt{2.16e-6} & \texttt{1.03e-3} & 7.439 \\
 &  & $10^{-7}$ & 140 & \texttt{9.82e-5} & \texttt{1.66e-6} & \texttt{3.40e-4} & 13.018 \\
\cmidrule(lr){2-8}
 & OGAProx & -- & 293 & \texttt{9.59e-5} & \texttt{2.20e-6} & \texttt{3.98e-4} & 0.739 \\
\midrule
$10^{-6}$ & AALM & $10^{-3}$ & 76 & \texttt{2.63e-7} & \texttt{2.73e-7} & \texttt{5.76e-5} & \textbf{0.323} \\
 &  & $10^{-5}$ & 51 & \texttt{1.27e-7} & \texttt{7.16e-7} & \texttt{2.86e-4} & 0.394 \\
 &  & $10^{-7}$ & \textbf{50} & \texttt{1.66e-7} & \texttt{7.76e-7} & \texttt{3.36e-4} & 0.619 \\
\cmidrule(lr){2-8}
 & ConEx & -- & 1643 & \texttt{3.92e-7} & \texttt{5.83e-7} & \texttt{3.42e-7} & 5.028 \\
\cmidrule(lr){2-8}
 & AGRPDA-L & $10^{-3}$ & 303 & \texttt{8.76e-7} & \texttt{8.01e-8} & \texttt{6.69e-5} & 5.602 \\
 &  & $10^{-5}$ & 275 & \texttt{8.96e-7} & \texttt{8.17e-8} & \texttt{2.09e-5} & 11.605 \\
 &  & $10^{-7}$ & 275 & \texttt{8.79e-7} & \texttt{7.96e-8} & \texttt{2.10e-5} & 22.076 \\
\cmidrule(lr){2-8}
 & OGAProx & -- & 617 & \texttt{9.43e-7} & \texttt{4.96e-8} & \texttt{2.26e-5} & 1.563 \\
\midrule
$10^{-8}$ & AALM & $10^{-3}$ & 172 & \texttt{1.44e-9} & \texttt{7.50e-9} & \texttt{2.20e-6} & 0.677 \\
 &  & $10^{-5}$ & 80 & \texttt{1.85e-10} & \texttt{6.09e-9} & \texttt{1.78e-6} & \textbf{0.516} \\
 &  & $10^{-7}$ & \textbf{75} & \texttt{9.17e-11} & \texttt{7.62e-9} & \texttt{4.10e-6} & 0.829 \\
\cmidrule(lr){2-8}
 & ConEx & -- & 3622 & \texttt{5.90e-10} & \texttt{9.37e-9} & \texttt{3.34e-9} & 10.582 \\
\cmidrule(lr){2-8}
 & AGRPDA-L & $10^{-3}$ & 579 & \texttt{7.34e-9} & \texttt{2.43e-9} & \texttt{3.88e-6} & 8.698 \\
 &  & $10^{-5}$ & 519 & \texttt{6.46e-9} & \texttt{2.74e-9} & \texttt{4.03e-6} & 17.153 \\
 &  & $10^{-7}$ & 508 & \texttt{5.92e-9} & \texttt{3.60e-9} & \texttt{1.13e-6} & 35.075 \\
\cmidrule(lr){2-8}
 & OGAProx & -- & 1211 & \texttt{8.48e-9} & \texttt{1.50e-9} & \texttt{1.22e-6} & 3.030 \\
\bottomrule\end{tabular}}
\begin{minipage}{\linewidth}\footnotesize
{``--'': no iterative inner solve.}
\end{minipage}
\end{table*}

\begin{table*}[t]\centering
\caption{Composite regression with $\mu=0$.}
\label{tab:comp2}
\resizebox{\textwidth}{!}{%
\begin{tabular}{clcrrrrr}\toprule
$\varepsilon_{\rm out}$ & Method & $e_k$ & $K$ & $|f(x_K)-f_{\rm ref}|$ & $\|[Cx_K-d]_+\|$ & $\|\Delta x_K\|+\|\Delta\lambda_K\|$ & Time (s) \\
\midrule
$10^{-4}$ & AALM & $10^{-3}$ & 34 & \texttt{6.20e-6} & \texttt{9.05e-5} & \texttt{1.06e-3} & \textbf{0.261} \\
 &  & $10^{-5}$ & \textbf{32} & \texttt{5.97e-6} & \texttt{8.79e-5} & \texttt{9.92e-4} & 0.266 \\
 &  & $10^{-7}$ & \textbf{32} & \texttt{5.93e-6} & \texttt{8.72e-5} & \texttt{9.89e-4} & 0.373 \\
\cmidrule(lr){2-8}
 & ConEx & -- & 10000$^{\dagger}$ & \texttt{4.18e-5} & \texttt{1.17e-3} & \texttt{1.42e-6} & 31.186 \\
\cmidrule(lr){2-8}
 & GRPDA-L & $10^{-3}$ & 122 & \texttt{2.46e-6} & \texttt{9.16e-5} & \texttt{2.26e-4} & 5.202 \\
 &  & $10^{-5}$ & 119 & \texttt{2.37e-6} & \texttt{9.73e-5} & \texttt{2.42e-4} & 7.304 \\
 &  & $10^{-7}$ & 119 & \texttt{2.36e-6} & \texttt{9.67e-5} & \texttt{2.42e-4} & 10.392 \\
\cmidrule(lr){2-8}
 & OGAProx & -- & 199 & \texttt{2.33e-6} & \texttt{9.57e-5} & \texttt{2.26e-4} & 0.788 \\
\midrule
$10^{-6}$ & AALM & $10^{-3}$ & \textbf{131} & \texttt{2.12e-8} & \texttt{9.54e-7} & \texttt{9.82e-6} & 0.733 \\
 &  & $10^{-5}$ & \textbf{131} & \texttt{4.73e-8} & \texttt{9.32e-7} & \texttt{1.14e-5} & \textbf{0.603} \\
 &  & $10^{-7}$ & 134 & \texttt{4.78e-8} & \texttt{9.35e-7} & \texttt{1.15e-5} & 0.935 \\
\cmidrule(lr){2-8}
 & ConEx & -- & 10000$^{\dagger}$ & \texttt{4.18e-5} & \texttt{1.17e-3} & \texttt{1.42e-6} & 31.186 \\
\cmidrule(lr){2-8}
 & GRPDA-L & $10^{-3}$ & 463 & \texttt{1.67e-8} & \texttt{9.81e-7} & \texttt{2.00e-6} & 10.869 \\
 &  & $10^{-5}$ & 419 & \texttt{4.27e-8} & \texttt{9.55e-7} & \texttt{2.49e-6} & 12.381 \\
 &  & $10^{-7}$ & 427 & \texttt{4.31e-8} & \texttt{9.46e-7} & \texttt{2.12e-6} & 19.391 \\
\cmidrule(lr){2-8}
 & OGAProx & -- & 713 & \texttt{8.31e-10} & \texttt{9.96e-7} & \texttt{1.32e-6} & 2.996 \\
\midrule
$10^{-8}$ & AALM & $10^{-3}$ & 310 & \texttt{3.48e-11} & \texttt{9.91e-9} & \texttt{1.66e-7} & 1.754 \\
 &  & $10^{-5}$ & \textbf{270} & \texttt{1.46e-10} & \texttt{9.81e-9} & \texttt{1.73e-7} & \textbf{1.072} \\
 &  & $10^{-7}$ & 278 & \texttt{1.52e-10} & \texttt{9.66e-9} & \texttt{1.72e-7} & 1.552 \\
\cmidrule(lr){2-8}
 & ConEx & -- & 10000$^{\dagger}$ & \texttt{4.18e-5} & \texttt{1.17e-3} & \texttt{1.42e-6} & 31.186 \\
\cmidrule(lr){2-8}
 & GRPDA-L & $10^{-3}$ & 1892 & \texttt{3.60e-12} & \texttt{9.97e-9} & \texttt{1.91e-8} & 34.467 \\
 &  & $10^{-5}$ & 1767 & \texttt{1.78e-10} & \texttt{9.79e-9} & \texttt{1.47e-8} & 31.517 \\
 &  & $10^{-7}$ & 1776 & \texttt{1.88e-10} & \texttt{9.79e-9} & \texttt{1.80e-8} & 37.632 \\
\cmidrule(lr){2-8}
 & OGAProx & -- & 2405 & \texttt{7.94e-11} & \texttt{9.89e-9} & \texttt{1.27e-8} & 9.023 \\
\bottomrule\end{tabular}}
\begin{minipage}{\linewidth}\footnotesize
  ``--'': no iterative inner solve. $\dagger$: outer cap reached; the terminal values are repeated.
\end{minipage}
\end{table*}

\subsection{Composite regression with affine inequalities}
We use a weighted elastic-net variant of the constrained lasso \cite[Eq.~(1)]{Gaines2018}, with an additional box constraint:
\[
\begin{gathered}
 \min_{x\in X} f(x):=\frac12\|Bx-y\|_2^2
       +\sum_{i=1}^n w_i|x_i|+\frac\mu2\|x\|_2^2
 \\ \text{subject to}\quad Cx-d\le0,\qquad X=[-a,a]^n.
\end{gathered}
\]
where $w_i>0$ are the weights.
In our notation, $\phi(x)=\frac12\|Bx-y\|_2^2$ and
$h(x)=\sum_iw_i|x_i|+\frac\mu2\|x\|_2^2+\iota_X(x)$.

We set $m=5000$, $n=3000$, and $q=40$, where
$B\in\R^{m\times n}$ and $C\in\R^{q\times n}$.
Both matrices have independent standard Gaussian entries before rescaling
to $\|B\|=\|C\|=1$. We take $a=1$ and
$w_i\sim\mathcal U[0.008,0.012]$, and examine $\mu=0.05$ and $\mu=0$.
For each value of $\mu$, we construct a synthetic instance with a known KKT pair. The matrices $B$, $C$, and the weights $w$ are shared between the two settings, while the KKT construction determines the corresponding $y$ and $d$. Within each setting, all methods use identical data and start
from $x_0=0$ and $\lambda_0=0$. The reference value $f_{\rm ref}$ is
accepted when its numerical Fenchel primal-dual objective-error certificate
is at most $10^{-12}$. This construction also applies at $\mu=0$.
 The two settings compare the presence and absence of explicit quadratic regularization. The Gaussian matrix $B$ has full column rank almost surely, so the least-squares term remains strongly convex at $\mu=0$.

We compare AALM with the deterministic constraint extrapolation method \cite[Algo.~1]{ConEx2023} (ConEx), the accelerated golden ratio primal-dual algorithm with linesearch \cite[Algo.~4.1]{Chang2022} (AGRPDA-L), and the optimistic gradient ascent--proximal point algorithm \cite[Sect.~1.2, Eqs.~(7), (8), (26)]{Bot2023} (OGAProx). We use their published outer-update formulas and parameter recurrences, with the primal proximal subproblems of the golden ratio methods solved inexactly as specified below.
For AALM, $(\alpha,\gamma,\beta,\sigma)=(100,98,0.1,3)$ and
$\tau=\min\{1,L_\phi^{-1}\}$, where $L_\phi=\|B\|^2$.
ConEx uses $\alpha_0=\mu$ and $M=2\|C\|_F$:
policy~(2.10) of~\cite{ConEx2023}, with $t_0=4L_\phi/\mu+2$, is used
when $\mu>0$; deterministic policy~(2.19), with $B_{\rm ConEx}=2$, is
used when $\mu=0$.
AGRPDA-L uses $\psi=1.5$, $\beta_0=\tau_0=1$, and linesearch reduction
$0.7$. At $\mu=0$, we use the golden ratio primal-dual algorithm with linesearch \cite[Algo.~3.1]{Chang2022} (GRPDA-L), with constant $\beta$ and linesearch safety factor $0.99$.
 OGAProx is applied to the negated Lagrangian with the variable roles exchanged: its coupling is $-\phi(x)-\langle\lambda,Cx-d\rangle+\iota_{\R_+^q}(\lambda)$ and its regularizer is $h(x)$. We use $\tau_0=1$, $\sigma_0=0.2$, $\theta_0=1$, $c_\alpha=2$, and $\nu=\mu$.
Its original recurrence yields constant step parameters at
$\mu=0$~\cite[Proposition~8]{Bot2023}.

The proximal map needed by ConEx and OGAProx is explicit:
\[ \operatorname{prox}_{t h}(v)
 =\Pi_X\!\left(\frac{\mathcal S_{tw}(v)}{1+t\mu}\right),
 \qquad
 [\mathcal S_{tw}(v)]_i=\operatorname{sign}(v_i)(|v_i|-tw_i)_+.
\]
Thus these two methods use weighted soft thresholding, quadratic scaling,
and coordinatewise clipping, without iterative inner solves.
AALM and the golden ratio methods use the same strongly convex version of
FISTA~\cite{Beck2017} with backtracking. For these solves, write
$F_k=S_k+\chi$, where
$\chi(x)=\sum_iw_i|x_i|+\iota_X(x)$ and $S_k$ contains all smooth terms,
including the quadratic regularizer and the proximal quadratic.
Their known strong-convexity lower bounds are $\mu+1/\tau$ for AALM and
$\mu+1/\tau^G_{k-1}$ for the golden ratio methods, respectively, where
$\tau^G_{k-1}$ is its current primal step. Both remain positive at $\mu=0$.
For each iterative inner solve, we use the fixed stopping conditions
\[
\begin{gathered}
 \operatorname{dist}(0,\partial F_k(x_{k+1}))\le e_k,\qquad
\text{at most }10^6\text{ inner updates}.
\end{gathered}
\]
The stationarity residual includes the normal cone of $X$.
As in sparse recovery, we use fixed inner tolerances for the finite-accuracy comparison and accept an inner solve when its stationarity test is satisfied. The outer stopping rule is
\begin{equation}\label{eq:outer-stop}
 R_K:=|f(x_K)-f_{\rm ref}|+\|[Cx_K-d]_+\|
 \le \varepsilon_{\rm out},\qquad
 \varepsilon_{\rm out}\in\{10^{-4},10^{-6},10^{-8}\}.
\end{equation}

Tables~\ref{tab:comp1}
and~\ref{tab:comp2} compare the cases $\mu=0.05$ and $\mu=0$,
respectively. Each successful row records the first outer iterate satisfying
\eqref{eq:outer-stop}: its iteration count $K$, absolute objective gap,
constraint violation, $\|\Delta x_K\|+\|\Delta\lambda_K\|$, and time,
where $\Delta x_K=x_K-x_{K-1}$ and
$\Delta\lambda_K=\lambda_K-\lambda_{K-1}$.
ConEx and OGAProx are independent of $e_k$ and are listed
only once for each outer threshold. The smallest successful $K$ and runtime are highlighted separately within each outer-tolerance group. An outer cap of $10^4$ is used for the reported ConEx run at $\mu=0$; the corresponding terminal values are marked when the stopping rule is not met.

Across Tables~\ref{tab:comp1} and~\ref{tab:comp2}, For each target, the best tested AALM settings give the smallest outer iteration count and the shortest reported runtime. Its advantage persists when the explicit quadratic regularizer is removed. ConEx reaches its outer cap at $\mu=0$ with a small successive difference but a larger feasibility residual, illustrating why the stopping rule monitors both objective accuracy and feasibility.  Tighter inner tolerances can reduce the outer count, especially for $\mu=0.05$ at high accuracy, but need not minimize runtime. For the strictest target, ${e_k}=10^{-5}$ gives the shortest AALM runtime in both settings. The comparisons therefore support the method's effectiveness with moderately accurate inner solves under the stated parameter choices.

\begin{table}[!ht]
\caption{QCQP comparison}
\label{tab:qcqp-ncg}
\resizebox{\linewidth}{!}{%
\begin{tabular}{clrrrrr}
\toprule
$\varepsilon_{\rm out}$ & Method & $K$ & $|f(x_K)-f_{\rm ref}|$ & $\|[g(x_K)]_+\|$ & $\|\Delta x_K\|+\|\Delta\lambda_K\|$ & Time (s) \\
\midrule
$10^{-4}$ & AALM & \textbf{37} & \texttt{7.07e-05} & \texttt{2.89e-05} & \texttt{6.72e-03} & \textbf{0.011} \\
\cmidrule(lr){2-7}
 & ACGD & 44 & \texttt{5.06e-05} & \texttt{4.79e-05} & \texttt{6.13e-03} & 0.014 \\
\cmidrule(lr){2-7}
 & APMM & 56 & \texttt{2.19e-05} & \texttt{7.72e-05} & \texttt{6.69e-03} & 0.012 \\
\midrule
$10^{-6}$ & AALM & \textbf{145} & \texttt{9.04e-07} & \texttt{8.61e-08} & \texttt{1.78e-03} & \textbf{0.031} \\
\cmidrule(lr){2-7}
 & ACGD & 262 & \texttt{7.66e-08} & \texttt{9.21e-07} & \texttt{8.00e-04} & 0.045 \\
\cmidrule(lr){2-7}
 & APMM & 541 & \texttt{2.36e-07} & \texttt{7.62e-07} & \texttt{3.06e-05} & 0.085 \\
\midrule
$10^{-8}$ & AALM & \textbf{273} & \texttt{2.96e-09} & \texttt{7.02e-09} & \texttt{2.35e-04} & \textbf{0.046} \\
\cmidrule(lr){2-7}
 & ACGD & 10000$^{\dagger}$ & \texttt{2.18e-08} & \texttt{1.06e-07} & \texttt{5.94e-06} & 1.450 \\
\cmidrule(lr){2-7}
 & APMM & 5396 & \texttt{2.37e-09} & \texttt{7.62e-09} & \texttt{5.13e-07} & 0.829 \\
\bottomrule
\end{tabular}}
\par\smallskip\footnotesize $\dagger$: target not reached.
\end{table}
 
\subsection{Quadratically constrained quadratic programming}
Following the quadratically constrained quadratic programming (QCQP) model in \cite{ZhangLanACGD}, we consider
\[ \min_{x\in X}\ f(x):=\tfrac12x^\top Qx+c^\top x,
 \qquad
 g_i(x):=\tfrac12x^\top H_ix+d_i^\top x\leq0,
 \quad i=1,\ldots,m,
\]
where $X=[-b,b]^n$, $Q\succeq0$, and $H_i\succ0$. We use one fixed instance with $n=100$, $m=10$,
$b=10$, and random seed $101$.
The matrices $Q$ and $H_i$ are randomly generated and scaled to
satisfy $\|Q\|=1$ and
$\|H_i\|\in[1,10]$.
The linear terms $c$ and $d_i$ are chosen to admit a known KKT pair
$(x^\star,\lambda^\star)$, where $x^\star\in\operatorname{int}X$,
$\|x^\star\|=1$, and $3$ quadratic constraints are active.
By convexity, $f_{\rm ref}=f(x^\star)$ is the optimal value.
All methods use the same instance, with all primal and dual states
initialized at zero.

We compare AALM (Algorithm \ref{alg:inexact}), ACGD in \cite{ZhangLanACGD}, and APMM in
\cite{DengInformsJoc}.
AALM uses $\phi=f$, $h=\iota_X$,
$(\alpha,\gamma,\beta,\sigma)=
(8,6,0.1,1)$, and
$\tau=1$, with $L_\phi=\|Q\|_2$ and $\tau L_\phi\leq1$.
ACGD uses the general convex parameter rule \cite[(2.16)]{ZhangLanACGD}, with $t_k=0.5(k-1)$, $\theta_k=(k-1)/k$, $\eta_k=2 L_C/k$ and $\omega_k=k$,
where
$L_C=\|Q+\sum_i\lambda_i^\star H_i\|_2+
10^{-3}\|(\|H_i\|_2)_{i=1}^m\|_2$.
This curvature bound uses the known reference multiplier.
APMM uses $\alpha_k=2/(k+1)$,
$X_k=X$, the known value $f_{\rm ref}$, and its prescribed output
selection rule. The AALM primal subproblem is solved by the projected Newton approach
for box constraints \cite{BertsekasSiam}. The linearly constrained quadratic subproblems of ACGD and APMM are
solved by a feasible active-set method. All three methods use the same absolute stationarity requirement for their respective original subproblem objectives:
$
 \operatorname{dist}(0,\partial F_k(x_{k+1}))\leq 10^{-5}.
$
We evaluate AALM at its generated primal iterates, ACGD at its weighted primal averages, and APMM at its prescribed primal outputs. With $x_k$ denoting the corresponding point, we record the first iteration satisfying
\[ R_K:=|f(x_K)-f_{\rm ref}|+\|[g(x_K)]_+\|
       \le \varepsilon_{\rm out},
 \qquad\varepsilon_{\rm out}\in\{10^{-4},10^{-6},10^{-8}\}.
\] 
 
 Table~\ref{tab:qcqp-ncg} shows that AALM reaches the reported targets in fewer outer iterations and less time than the successful comparison runs. At $10^{-8}$, ACGD reaches the iteration cap, and APMM takes substantially more outer iterations than AALM.

\section{Conclusions}\label{sec:conclusion}
In this paper, we have developed a class of augmented Lagrangian methods with Nesterov acceleration for composite convex optimization with nonlinear inequality constraints. Under weighted summability of the primal stationarity errors, the methods achieve simultaneous $\mathcal O(k^{-2})$ rates for the objective residual and constraint violation at the generated primal iterates, including the case $\alpha=3$. For $\alpha>3$ and $2<\gamma<\alpha-1$, both residuals are $o(k^{-2})$, the entire primal-dual sequence converges to a KKT pair, and successive primal and dual differences are $o(k^{-1})$. Compared with general convex augmented Lagrangian methods \cite{XuLinear2021,XuMathProg} and existing accelerated methods for nonlinear convex inequalities \cite{HamedaniAybat2021,ConEx2023,ZhuSiopt,LinNeurips,ZhangLanACGD,DengInformsJoc}, the proposed methods combine accelerated last-iterate rates with full primal-dual convergence while allowing inexact primal solves and requiring neither the objective nor the constraints to be strongly convex. The convergence analysis extends related results on accelerated decay and iterate convergence for unconstrained problems \cite{AttouchPtSiopt,ChambolleJota} and discrete augmented Lagrangian methods for affine equality constraints \cite{HeArxiv} to nonlinear convex inequalities. Numerical experiments demonstrate the computational advantages of the proposed methods on the tested problems.

\appendix

\makeatletter
\let\aalm@original@seccntformat\@seccntformat
\renewcommand{\@seccntformat}[1]{\protect{\appendixname~\csname the#1\endcsname:}\enspace}
\makeatother
\section[Technical lemmas]{Technical lemmas}\label{sec:appendix}

In this part, we collect the technical lemmas used in the convergence analysis. 

The first is the standard variational characterization of metric projections \cite[Theorem~3.16]{Bauschke16}.

\begin{lemma}\label{lem:proj}
Let $C\subseteq\R^m$ be a nonempty closed convex set. For every
$z,y\in\R^m$,
\[
  y=\proj_C(z)
  \quad\Longleftrightarrow\quad
  y\in C,\qquad
  \ip{z-y}{v-y}\leq0
  \quad\forall\,v\in C.
\]
Moreover, for every $\theta>0$,
\[
  y=\proj_C
   \bigl(\theta z+(1-\theta)y\bigr)
  \quad\Longleftrightarrow\quad
  y=\proj_C(z).
\]
\end{lemma}
 
\begin{lemma}[{Opial's lemma \cite{Bauschke16}}]\label{lem:opial_dis}
Let $S\subset\R^d$ be nonempty and let $\{u_k\}\subset\R^d$.
Fix any norm on $\R^d$. Suppose that $\lim_{k\to+\infty}\|u_k-v\|$
exists and is finite for every $v\in S$, and that every cluster point of
$\{u_k\}$ belongs to $S$. Then $u_k$ converges to a point of $S$.
\end{lemma}

\begin{lemma}\cite[Lemma~5.14]{AttouchMathProg}\label{le_disc_per}
Let $\{a_k\}_{k\geq1}$ and $\{b_k\}_{k\geq1}$ be two nonnegative sequences.
Assume that $\sum_{k=1}^{+\infty}b_k<+\infty$ and
$
    a_k^2\leq c^2+\sum_{j=1}^{k}b_ja_j
$ for every $k\geq1$,
where $c\geq0$. Then
$
    \sup_{k\geq1}a_k
    \leq c+\sum_{j=1}^{+\infty}b_j
    <+\infty.
$
\end{lemma}
 
\begin{lemma}\label{lem:bounded_dis}
Let $\{x_k\}_{k\geq0}\subset\R^n$, $x^\star\in\R^n$, and $\gamma>0$. Suppose that
\[
  \sup_{k\geq0}
  \norm{
    \gamma(x_{k+1}-x^\star)
    +k(x_{k+1}-x_k)
  }
  <+\infty.
\]
Then $\{x_k\}$ is bounded and
$
  \norm{x_{k+1}-x_k}
  =
  \mathcal O(k^{-1}).
$
\end{lemma}
\begin{proof}
By assumption, there exists $M>0$ such that
\[
  \norm{
    \gamma(x_{k+1}-x^\star)
    +k(x_{k+1}-x_k)
  }
  \leq M
\]
for all $k\geq0$. Rearranging gives
\[
  (k+\gamma)(x_{k+1}-x^\star)
  =
  k(x_k-x^\star)
  +
  \gamma(x_{k+1}-x^\star)
  +k(x_{k+1}-x_k).
\]
Hence
\[
  \norm{x_{k+1}-x^\star}
  \leq
  \frac{k}{k+\gamma}\norm{x_k-x^\star}
  +\frac{M}{k+\gamma}.
\]
Let
$
  C:=\max\{\norm{x_0-x^\star},M/\gamma\}.
$ We show by induction that
$
  \norm{x_k-x^\star}\leq C
$
for all $k\geq0$. The claim is clear for $k=0$. If
$
  \norm{x_k-x^\star}\leq C,
$
then, since $M\leq\gamma C$,
\[
  \norm{x_{k+1}-x^\star}
  \leq
  \frac{kC+M}{k+\gamma}
  \leq
  \frac{kC+\gamma C}{k+\gamma}
  =C.
\]
Therefore, $\{x_k\}$ is bounded. Finally, for $k\geq1$,
\[
  k\norm{x_{k+1}-x_k}
  \leq
  M+\gamma\norm{x_{k+1}-x^\star}
  \leq
  M+\gamma C.
\]
Therefore,
$
  \norm{x_{k+1}-x_k}
  =
  \mathcal O(k^{-1}),
$
which completes the proof.
\end{proof}

The following lemma controls the positive part of a sequence satisfying a discrete summation identity.

\begin{lemma}\label{lem:dis_vol}
Let $\{G_k\}_{k\geq1}$, $\{B_k\}_{k\geq1}$, and
$\{M_k\}_{k\geq1}$ be sequences in $\R^m$, and let
$\{a_k\}_{k\geq2}$ be nonnegative.
Suppose that $\{M_k\}$ is componentwise nondecreasing,
$M_1=0$, and
\[
  G_k+\sum_{j=2}^k a_jG_j=B_k-M_k,
  \qquad k\geq1.
\]
Then the following statements hold.
\begin{itemize}
  \item[(i)] If $\sup_{k\geq1}\|B_k\|\leq C$, then
  \[
    \sup_{k\geq1}\|\pos{G_k}\|\leq2C.
  \]

  \item[(ii)] If $\{B_k\}$ has a finite limit and
  $\sum_{k=2}^{+\infty}a_k=+\infty$, then
  \[
    \lim_{k\to+\infty}\|\pos{G_k}\|=0.
  \]
  Moreover, if for some component $i$ there exists an
  integer $K_i\geq1$ such that
  $M_{k,i}=M_{K_i,i}$ for all $k\geq K_i$, then
  \[
    \lim_{k\to+\infty}G_{k,i}=0.
  \]
\end{itemize}
\end{lemma}

\begin{proof}
Set
\[
  Z_k:=\sum_{j=2}^k a_jG_j,
  \qquad
  \theta_k:=\frac1{1+a_k}\in(0,1].
\]
Then $Z_1=0$, and the assumed relation gives
\begin{equation}\label{eq_GBZ}
	G_k=B_k-M_k-Z_k.
\end{equation}
For $k\geq2$,
\[
  Z_k
  =Z_{k-1}+a_kG_k
  =Z_{k-1}+a_k(B_k-M_k-Z_k).
\]
Rearranging gives
\[
  Z_k
  =\theta_kZ_{k-1}+(1-\theta_k)(B_k-M_k).
\]
Iterating this identity and using $Z_1=0$, we obtain
\begin{equation}\label{eq_ZKJ}
  Z_k=\sum_{j=2}^k\omega_{j,k}(B_j-M_j),
  \qquad\text{with} \qquad
  \omega_{j,k}:=
  (1-\theta_j)\prod_{\ell=j+1}^k\theta_\ell,
  \quad 2\leq j\leq k.
\end{equation}
The weights are nonnegative and satisfy
\begin{equation}\label{eq_wjk}
  \begin{aligned}
    \sum_{j=2}^k\omega_{j,k}
=
    \sum_{j=2}^k
    \left(
      \prod_{\ell=j+1}^k\theta_\ell
      -\prod_{\ell=j}^k\theta_\ell
    \right)=1-\prod_{\ell=2}^k\theta_\ell\in[0,1].
  \end{aligned}
\end{equation}
Consequently, \eqref{eq_GBZ} and \eqref{eq_ZKJ} give
\begin{equation}\label{eq:dis_vol_rep}
  G_k
  =
  B_k-\sum_{j=2}^k\omega_{j,k}B_j
  -M_k+\sum_{j=2}^k\omega_{j,k}M_j.
\end{equation}

(i)
Since $M_1=0$ and $\{M_k\}$ is componentwise
nondecreasing, $0\leq M_j\leq M_k$ for $j\leq k$.
Thus \eqref{eq_wjk} gives, componentwise,
\[
  -M_k+\sum_{j=2}^k\omega_{j,k}M_j
  \leq
  -M_k+\left(\sum_{j=2}^k\omega_{j,k}\right)M_k
  =
  -\left(\prod_{\ell=2}^k\theta_\ell\right)M_k
  \leq0.
\]
Combining this inequality with \eqref{eq:dis_vol_rep}
and taking positive parts yields
\begin{equation}\label{eq:dis_vol_bound}
  \|\pos{G_k}\|
  \leq
  \left\|\pos{B_k-\sum_{j=2}^k\omega_{j,k}B_j}\right\|
  \leq
  \left\|B_k-\sum_{j=2}^k\omega_{j,k}B_j\right\|.
\end{equation}
If $\sup_{k\geq1}\|B_k\|\leq C$, then
\[
  \|\pos{G_k}\|
  \leq
  \|B_k\|+\sum_{j=2}^k\omega_{j,k}\|B_j\|
  \leq
  C\left(1+\sum_{j=2}^k\omega_{j,k}\right)
  \leq2C,
\]
where the last inequality follows from \eqref{eq_wjk}. Taking the supremum over $k\geq1$ proves~(i).

(ii) Since $a_\ell\geq0$, induction gives
\[
  \prod_{\ell=2}^k(1+a_\ell)
  \geq1+\sum_{\ell=2}^k a_\ell,
  \qquad k\geq2.
\]
Indeed, equality holds for $k=2$. If the inequality
holds for $k-1$, then
\[
  \prod_{\ell=2}^k(1+a_\ell)
  \geq
  \left(1+\sum_{\ell=2}^{k-1}a_\ell\right)(1+a_k)
  \geq
  1+\sum_{\ell=2}^k a_\ell.
\]
Since $\theta_\ell=1/(1+a_\ell)$ and
$\sum_{\ell=2}^{+\infty}a_\ell=+\infty$, it follows that
\begin{equation}\label{eq_prodTheta}
  0<
  \prod_{\ell=2}^k\theta_\ell
  =
  \frac{1}{\prod_{\ell=2}^k(1+a_\ell)}
  \leq
  \frac{1}{1+\sum_{\ell=2}^k a_\ell}
  \to0, \qquad\text{as }k\to+\infty.
\end{equation}
For each fixed $j\geq2$, we also have
$\sum_{\ell=j+1}^{+\infty}a_\ell=+\infty$.
Applying the same estimate to the product starting
at $j+1$, from \eqref{eq_ZKJ}, we obtain
\begin{equation}\label{eq_wjkk}
  0\leq\omega_{j,k}
  =
  (1-\theta_j)\prod_{\ell=j+1}^k\theta_\ell
  \leq
  \frac{1-\theta_j}
       {1+\sum_{\ell=j+1}^k a_\ell}
  \to0
  \qquad\text{as }k\to+\infty.
\end{equation}

Since $\{B_k\}$ has a finite limit, let $B^\infty=\lim_{k\to+\infty}B_k$. By
\eqref{eq_wjk},
\[
  \sum_{j=2}^k\omega_{j,k}B_j-B^\infty
  =
  \sum_{j=2}^k\omega_{j,k}(B_j-B^\infty)
  -
  \left(\prod_{\ell=2}^k\theta_\ell\right)B^\infty.
\]
For every $\varepsilon>0$, choose $N\geq2$ such that
$\|B_j-B^\infty\|\leq\varepsilon$ for all $j>N$.
It follows from \eqref{eq_wjk} that for any  $k>N$,
\begin{align*}
  \left\|
    \sum_{j=2}^k\omega_{j,k}B_j-B^\infty
  \right\|
  &\leq
  \sum_{j=2}^N\omega_{j,k}\|B_j-B^\infty\|+
  \varepsilon\sum_{j=N+1}^k\omega_{j,k}
  +\left(\prod_{\ell=2}^k\theta_\ell\right)\|B^\infty\|\\
  &\leq
  \sum_{j=2}^N\omega_{j,k}\|B_j-B^\infty\|
  +\varepsilon
  +\left(\prod_{\ell=2}^k\theta_\ell\right)\|B^\infty\|.
\end{align*}
For this fixed $N$, \eqref{eq_wjkk} gives
$\lim_{k\to+\infty}\omega_{j,k} = 0$ for any $2\leq j\leq N$. Then from \eqref{eq_prodTheta},
\[
  0\leq
  \limsup_{k\to+\infty}
  \left\|
    \sum_{j=2}^k\omega_{j,k}B_j-B^\infty
  \right\|
  \leq\varepsilon.
\]
As $\varepsilon>0$ is arbitrary, this upper limit is
zero. Therefore,
\begin{equation}\label{eq:dis_weighted_B}
  \lim_{k\to+\infty}
  \sum_{j=2}^k\omega_{j,k}B_j=B^\infty.
\end{equation}
Using \eqref{eq:dis_vol_bound} and $B^\infty=\lim_{k\to+\infty}B_k$, we then obtain
\[
  0\leq\|\pos{G_k}\|
  \leq
  \|B_k-B^\infty\|
  +
  \left\|
    \sum_{j=2}^k\omega_{j,k}B_j-B^\infty
  \right\|
  \to0,\quad \text{as}\quad k\to+\infty.
\]

Finally, suppose that $M_{k,i}=M_{K_i,i}$ for every
$k\geq K_i$. By \eqref{eq_wjk}, for $k\geq K_i$,
\begin{align*}
  \sum_{j=2}^k\omega_{j,k}M_{j,i}-M_{K_i,i}
  &=
  \sum_{j=2}^k\omega_{j,k}
  (M_{j,i}-M_{K_i,i})
  -\left(\prod_{\ell=2}^k\theta_\ell\right)M_{K_i,i}\\
  &=
  \sum_{j=2}^{K_i-1}\omega_{j,k}
  (M_{j,i}-M_{K_i,i})
  -\left(\prod_{\ell=2}^k\theta_\ell\right)M_{K_i,i}.
\end{align*}
The finite sum tends to zero because
$\omega_{j,k}\to0$ for each fixed $j$ (by \eqref{eq_wjkk}), and the last
term tends to zero by \eqref{eq_prodTheta}.
Thus,
\begin{equation}\label{eq:limM}
  \lim_{k\to+\infty}
  \sum_{j=2}^k\omega_{j,k}M_{j,i}
  =M_{K_i,i}.
\end{equation}
For $k\geq K_i$, the $i$th component of
\eqref{eq:dis_vol_rep} can be written as
\[
  G_{k,i}
  =
  \left(
    B_{k,i}-\sum_{j=2}^k\omega_{j,k}B_{j,i}
  \right)
  +
  \left(
    \sum_{j=2}^k\omega_{j,k}M_{j,i}-M_{K_i,i}
  \right).
\]
Both terms on the right tend to zero by
\eqref{eq:dis_weighted_B} and  \eqref{eq:limM}.
Hence $\lim_{k\to+\infty}G_{k,i}=0$.
\end{proof}

\begin{lemma}\label{lem:est_diff}
Let $\{u_k\}_{k\geq0}\subset\R^{d_1}$ and
  $\{r_k\}_{k\geq1}\subset\R^{d_2}$ satisfy
  \[
    \sum_{k=1}^{+\infty}
    k\left(
      \|u_{k+1}-u_k\|^2+\|r_{k+1}\|^2
    \right)<+\infty.
  \]
If $\{n_j\}$ is a sequence of positive integers with
$n_j\to+\infty$ and $\rho>1$, then
  \[
    \sup_{n_j\leq q\leq\lfloor\rho n_j\rfloor}
    \|u_q-u_{n_j}\|\to0, \quad \text{as\ } n_j\to +\infty
  \]
  Moreover, there exist integers
  $q_j\in[n_j,\lfloor\rho n_j\rfloor]$ such that
  \[
    q_j\left(
      \|u_{q_j}-u_{q_j-1}\|
      +\|u_{q_j+1}-u_{q_j}\|
      +\|r_{q_j+1}\|
    \right)\to0, \quad \text{as\ } n_j\to +\infty
  \]
\end{lemma}
\begin{proof}
Set $N_j=\lfloor\rho n_j\rfloor$. For every integer
$q\in[n_j,N_j]$,   the Cauchy-Schwarz
inequality yields
\begin{align*}
  \|u_q-u_{n_j}\|
  &\leq
  \sum_{i=n_j+1}^{q}\|u_i-u_{i-1}\|\leq
  \left(
    \sum_{i=n_j+1}^{q}
    i\|u_i-u_{i-1}\|^2
  \right)^{1/2}
  \left(
    \sum_{i=n_j+1}^{q}\frac1i
  \right)^{1/2}.
\end{align*}
Moreover,
 \[
  \sum_{i=n_j+1}^{N_j}\frac1i
  \leq
  \int_{n_j}^{N_j}\frac{dt}{t}
  =
  \log\frac{N_j}{n_j}
  \leq\log\rho.
\]
Therefore,
\[
  \sup_{\substack{n_j\leq q\leq N_j}}
  \|u_q-u_{n_j}\|
  \leq
  \sqrt{\log\rho}
  \left(
    \sum_{i=n_j+1}^{+\infty}
    i\|u_i-u_{i-1}\|^2
  \right)^{1/2}.
\]
Since $
    \sum_{k=1}^{+\infty}
    k 
      \|u_{k+1}-u_k\|^2 <+\infty$, as $n_j\to+\infty$, 
its tail $\sum_{i=n_j+1}^{+\infty}
    i\|u_i-u_{i-1}\|^2$ tends to zero, proving the first assertion.

For the second assertion, define
\begin{equation}\label{eq_eta}
  \eta_j:=
  \sum_{q=n_j}^{N_j}q
  \left(
    \|u_q-u_{q-1}\|^2
    +\|u_{q+1}-u_q\|^2
    +\|r_{q+1}\|^2
  \right).
\end{equation}
By assumption, $\eta_j\geq0$, and all three weighted series appearing in the definition
of $\eta_j$ converge. Since $n_j\to+\infty$, their
sums over $q=n_j,\ldots,N_j$ tend to zero. Hence
\begin{equation}\label{eq_etaj}
	\lim_{j\to+\infty}\eta_j=0.
\end{equation}

Since $\{n_j,n_j+1,\ldots,N_j\}$ is finite and nonempty,
there exists an index $q_j$ in this set such that
\begin{align*}
  &q_j\left(
    \|u_{q_j}-u_{q_j-1}\|^2
    +\|u_{q_j+1}-u_{q_j}\|^2
    +\|r_{q_j+1}\|^2
  \right)\\
  &\qquad\leq
  q\left(
    \|u_q-u_{q-1}\|^2
    +\|u_{q+1}-u_q\|^2
    +\|r_{q+1}\|^2
  \right)
\end{align*}
for every $q\in\{n_j,n_j+1,\ldots,N_j\}$.
Summing this inequality over $q=n_j,\ldots,N_j$ and using
\eqref{eq_eta}, we obtain
\begin{equation}\label{eq_qjeq}
  q_j\left(
    \|u_{q_j}-u_{q_j-1}\|^2
    +\|u_{q_j+1}-u_{q_j}\|^2
    +\|r_{q_j+1}\|^2
  \right)\leq
  \frac{\eta_j}{N_j-n_j+1}
  \leq
  \frac{\eta_j}{(\rho-1)n_j},
\end{equation}
where the last inequality follows from
$N_j-n_j+1=\lfloor\rho n_j\rfloor-n_j+1
\geq(\rho-1)n_j$.

Multiplying \eqref{eq_qjeq} by $q_j$ and using $q_j\leq N_j\leq \rho n_j$ gives
\[q_j^2\left(
    \|u_{q_j}-u_{q_j-1}\|^2
    +\|u_{q_j+1}-u_{q_j}\|^2
    +\|r_{q_j+1}\|^2
  \right)\leq
  \frac{\rho \eta_j}{\rho-1}.
\]
Finally, applying
$(a+b+c)^2\leq3(a^2+b^2+c^2)$ and using \eqref{eq_etaj} yields
\begin{align*}
  0
 \leq
  \left[
    q_j\left(
      \|u_{q_j}-u_{q_j-1}\|
      +\|u_{q_j+1}-u_{q_j}\|
      +\|r_{q_j+1}\|
    \right)
  \right]^2\leq
  \frac{3\rho \eta_j}{\rho-1}
  \to0,
\end{align*}
as $j \to+\infty$. The required conclusion follows. This completes the proof.
\end{proof}

\begin{lemma}\label{lem:w_limit}
Let $\{h_k\}_{k\geq1}$ and $\{E_k\}_{k\geq1}$ be real
sequences, and let $\{R_k\}_{k\geq1}$ be nonnegative.
Suppose that, for some $a\ge 1$,
\[
  (k-1)(h_k-h_{k-1})+ah_k=E_k-R_k,
  \qquad k\geq2.
\]
If $\lim_{k\to+\infty}E_k$ exists and is finite, and
$\sum_{k=2}^{+\infty}R_k/k<+\infty$, then
\[
  \lim_{k\to+\infty}h_k
  =\frac1a\lim_{k\to+\infty}E_k.
\]
\end{lemma}

\begin{proof}
Set
\[
  P_1=1,\qquad
  P_k=\prod_{j=2}^k\left(1+\frac{a}{j-1}\right),
  \qquad k\geq2.
\]
By definition,
\begin{equation}\label{eq:weighted_product_identity}
  P_k=\frac{k+a-1}{k-1}P_{k-1},
  \qquad
  P_k-P_{k-1}
  =\frac{aP_{k-1}}{k-1}
  =\frac{aP_k}{k+a-1}>0.
\end{equation}
Thus $\{P_k\}$ is strictly increasing and $P_k\geq P_1=1$.
In particular,
\[
  P_k-P_{k-1}\geq\frac{a}{k-1}.
\]
Summing this inequality gives
\[
  P_k
  =1+\sum_{j=2}^k(P_j-P_{j-1})
  \geq1+a\sum_{j=2}^k\frac1{j-1}.
\]
This implies $P_k\to+\infty$.

We next derive a representation of $h_k$.
For every $j\geq2$, the assumed recurrence can be written as
\[
  (j+a-1)h_j=(j-1)h_{j-1}+E_j-R_j.
\]
Multiplying by $P_j/(j+a-1)$ and using \eqref{eq:weighted_product_identity}
we obtain
\[
  P_jh_j-P_{j-1}h_{j-1}
  =
  \frac{P_j-P_{j-1}}a E_j
  -\frac{P_jR_j}{j+a-1}.
\]
Summing from $j=2$ to $k$, using $P_1=1$, and dividing
by $P_k$, we obtain
\begin{equation}\label{eq:weighted_distance_rep}
  h_k
  =
  \frac{h_1}{P_k}
  +\frac1{aP_k}\sum_{j=2}^k(P_j-P_{j-1})E_j
  -\frac1{P_k}\sum_{j=2}^k\frac{P_jR_j}{j+a-1}.
\end{equation}

Since $a\ge 1$, $P_k$ is strictly increasing with
$P_k\to+\infty$, and $\{E_k\}$ has a finite limit,
we may apply the Stolz-Ces\`aro theorem
\cite{Choudary} to
\[
  a_k=\sum_{j=2}^k(P_j-P_{j-1})E_j,
  \qquad
  b_k=aP_k,
\]
to obtain
\begin{equation}\label{eq_term1}
\begin{aligned}
\lim_{k\to+\infty} \frac{\sum_{j=2}^k(P_j-P_{j-1})E_j}{aP_k}
  &=
  \lim_{k\to+\infty}
  \frac{a_k-a_{k-1}}{b_k-b_{k-1}}
  \\
  &=\lim_{k\to+\infty}
  \frac{(P_k-P_{k-1})E_k}{a(P_k-P_{k-1})}
  =
  \frac1a\lim_{k\to+\infty}E_k.
\end{aligned}
\end{equation}

For the term containing $R_j$, since $a\geq1$ and
$P_j\leq P_k$ for $j\leq k$, we have
\[
  0\leq
  \frac{P_jR_j}{P_k(j+a-1)}
  \leq\frac{R_j}{j},
  \qquad 2\leq j\leq k.
\]
Thus, for every fixed $N\geq2$ and $k>N$,
\begin{align*}
  0
  &\leq
  \frac1{P_k}\sum_{j=2}^k\frac{P_jR_j}{j+a-1}\leq
  \frac1{P_k}\sum_{j=2}^N\frac{P_jR_j}{j+a-1}
  +\sum_{j=N+1}^{+\infty}\frac{R_j}{j}.
\end{align*}
For fixed $N$, the first term on the right tends to zero
as $P_k\to+\infty$. Hence
\[
  0\leq
  \limsup_{k\to+\infty}
  \frac1{P_k}\sum_{j=2}^k\frac{P_jR_j}{j+a-1}
  \leq
  \sum_{j=N+1}^{+\infty}\frac{R_j}{j}.
\]
Letting $N\to+\infty$ and using the assumption
$\sum_{k=2}^{+\infty}R_k/k<+\infty$, we obtain
\[
  \lim_{k\to+\infty}
  \frac1{P_k}\sum_{j=2}^k\frac{P_jR_j}{j+a-1}=0.
\]
Together with $h_1/P_k\to0$,
\eqref{eq:weighted_distance_rep} and \eqref{eq_term1} yield
\[
  \lim_{k\to+\infty}h_k
  =\frac1a\lim_{k\to+\infty}E_k.
\]

\end{proof}

\end{document}